\documentclass[12pt]{article}
\usepackage{amsmath}
\newcommand{\bb}{\begin{eqnarray*}}

\newcommand{\E}{{\bf {E}}}
\newcommand{\p}{ { \bf P} }
\newcommand{\Z}{{\bf {Z}}}
\newcommand{\R}{{\bf {R}}}
\newcommand{\N}{{\mathbb{N}}}

\newcommand{\F}{{\cal F}}

\newcommand{\var}{{\rm Var}}
\newcommand{\mathbb}{\bf}

\newcommand{\ee}{\end{eqnarray*}}

\newtheorem{theorem}{Theorem}

\newtheorem{corollary}[theorem]{Corollary}

\newtheorem{definition}[theorem]{Definition}

\newtheorem{lemma}[theorem]{Lemma}

\newtheorem{proposition}[theorem]{Proposition}
\newtheorem{remark}[theorem]{Remark}
\newtheorem{comment}[theorem]{Comment}

\begin{document}
\begin{center}
{\bf \Large Moderate deviations for stationary sequences of bounded
random variables}\vskip15pt

J\'er\^ome Dedecker $^{a}$, Florence Merlev\`{e}de $^{b}$, Magda
Peligrad $^{c}$ \footnote{Supported in part by a Charles Phelps
Taft Memorial Fund grant and NSA\ grant, H98230-07-1-0016.} {\it
and\/} Sergey Utev   $^{d}$
\end{center}
\vskip10pt
$^a$ Universit\'e Paris 6, LSTA, 175 rue
du Chevaleret, 75013 Paris, FRANCE\\ \\
$^b$ Universit\'e Paris 6, LPMA and C.N.R.S UMR 7599, 175 rue
du Chevaleret, 75013 Paris, FRANCE\\ \\
$^c$ Department of Mathematical Sciences, University of Cincinnati,
PO Box 210025, Cincinnati, Oh 45221-0025,\\ \\
$^d$ School of Mathematical Sciences, University of Nottingham,
Nottingham, NG7 2RD, UK
\vskip10pt


\begin{center}
{\bf Abstract}\vskip10pt
\end{center}

In this paper we derive the moderate deviation principle for
stationary sequences of bounded random variables under
martingale-type conditions. Applications to functions of
$\phi$-mixing sequences, contracting Markov chains, expanding maps
of the interval, and symmetric random walks on the circle are given.
\section{Introduction}
$\quad \;$
For the stationary sequence $(X_i)_{ i \in {\mathbb Z}}$ of centered random variables,
define the partial sums and the normalized partial
sums process by
\begin{equation*}
S_{n}=\sum_{j=1}^{n}X_{j}\;\;\mbox{and}\;\;
W_{n}(t)=n^{-1/2}\sum_{i=1}^{[nt]}X_{i}\,.
\end{equation*}

In this paper we are concerned with the moderate deviation principle
for the normalized partial sums process $W_n$, considered as an
element of $D([0,1])$ (functions on $[0,1]$ with left-hand limits
and continuous from the right), equipped with the Skorohod topology
(see Section 14 in Billingsley (1968) for the description of the
topology on $D([0,1])$). More exactly, we say that the family of
random variables $\{W_n,n> 0\}$ satisfies the Moderate Deviation
Principle (MDP) in $D[0,1]$ with speed $a_n \rightarrow 0$ and good
rate function $I(.)$, if the level sets $\{ x , I(x) \leq \alpha \}$
are compact for all $\alpha < \infty$, and for all Borel sets
\begin{eqnarray}
-\inf_{t\in \Gamma^{0}}I(t)  &\leq&\lim\inf_{n}a_{n}\log {\bf P}(\sqrt{a_n}W_n\in \Gamma)
\nonumber\\
&\leq&\lim\sup_{n}a_{n}\log {\bf P}(\sqrt{a_n}W_n\in
\Gamma)\leq-\inf_{t\in\bar{\Gamma}}I(t) \, .\label{mdpkey}
\end{eqnarray}

The Moderate Deviation Principle is an intermediate behavior
between the central limit theorem ($a_n=a$) and Large Deviation ($a_n=a/n$).
Usually, MDP has a simpler rate function, inherited from the approximated
Gaussian process, and holds for a larger class of dependent
random variables than the large deviation principle.

De Acosta and Chen (1998) used the renewal theory to derive the MDP
for bounded functionals of geometrically ergodic stationary Markov
chains. Puhalskii (1994) and Dembo (1996) applied the stochastic
exponential to prove the MDP for martingales.  Starting from the
martingale case and using the so-called coboundary decomposition due
to Gordin (1969) ($X_k=M_k+Z_k-Z_{k+1}$, where $M_k$ is a stationary
martingale difference), Gao (1996) and Djellout (2002) obtained the
MDP for $\phi$-mixing sequences with summable mixing rate. In the
context of Markov chains, the coboundary decomposition is known as
the Poisson equation. Starting from this equation, Delyon, Juditsky
and Liptser (2006) proved the MDP for $n^{-1/2}
\sum_{k=1}^{n}H(Y_k)$, where $H$ is a Lipschitz function, and
$Y_n=F(Y_{n-1}, \xi_n)$, where $F$ satisfies $|F(x,z)-F(y,t)|\leq
\kappa |x-y|+ L|z-t|$ with $\kappa < 1$, and $(\xi_n)_{n \geq 1}$ is
an iid sequence of random variables independent of $Y_0$. In their
paper, the random variables are not assumed to   be bounded: the
authors only assume that there exists a positive $\delta$ such that
${\mathbf E} (e^{\delta |\xi_1|}) < \infty$. They strongly used the
Markov structure to derive some appropriate properties for the
coboundary (see their lemma 4.2).

In this paper we propose a modification of the martingale
approximation approach that allows to avoid the coboundary
decomposition and thus to enlarge the class of dependent sequences
known to satisfy the moderate deviation principle. Recent or new
exponential inequalities are applied to justify the martingale
approximation. The conditions involved in our results are well
adapted to a large variety of examples, including regular
functionals of linear processes, expanding maps of the interval
and symmetric random walks on the circle.

The paper is organized as follows. In Section 2
we state the main results. A discussion of the conditions,
clarifications,  and some simple examples and extensions
follow. Section 3 describes the applications, while Section
4 is dedicated to the proofs. Several technical lemmas are
proved in the appendix.

\section{Results}

$\quad \;$ From now on, we assume that the stationary sequence
$(X_i)_{ i \in {\mathbb Z}}$ is given by $X_i=X_0\circ T^i$, where
$T:\Omega \mapsto \Omega$ is a bijective bimeasurable transformation
preserving the probability $\p$ on $(\Omega,{\cal A})$. For a
subfield ${\cal F}_0 $ satisfying ${\cal F}_0 \subseteq T^{-1
}({\cal F}_0)$, let ${\mathcal F}_i=T^{-i}({\mathcal F}_0)$. By
$\Vert X\Vert _{\infty }$  we denote the $\mathbf L_\infty$-norm,
that is the smallest $u$ such that $\p (|X|>u)=0$.

Our first theorem and its corollary treat the so-called adapted case, $X_0$ being
${\mathcal F}_0$--measurable and so the sequence
$(X_i)_{ i \in {\mathbb Z}}$ is adapted to the
filtration $({\mathcal F}_i)_{ i \in {\mathbb Z}}$.
\begin{theorem}
\label{mainresult} Assume that $\|X_0\|_\infty< \infty$ and that
$X_0$ is ${\mathcal F}_0$--measurable. In addition, assume that
\begin{equation}
\sum_{n=1}^{\infty }n^{-3/2}\Vert \mathbf{E}(S_{n}|\mathcal{F}_{0})\Vert
_{\infty }<\infty \,,  \label{mw}
\end{equation}
and that there exists $\sigma ^{2}\geq 0$ with
\begin{equation}
\lim_{n\rightarrow \infty }\Vert
n^{-1}\mathbf{E}(S_{n}^{2}|\F_0)-\sigma ^{2}\Vert _{\infty}=0 \, .
\label{s2inf}
\end{equation}
Then, for all positive sequences $a_n$ with $a_n\to0$ and $na_n\to \infty$,
the normalized partial sums processes $W_{n}(.)$ satisfy (\ref{mdpkey})
with the good rate function $I_\sigma(\cdot)$ defined by
\begin{equation}
I_\sigma(h) =
\frac{1}{2\sigma^2}\int_{0}^{1}(h^{\prime}(u))^2du \quad
\label{rate}
\end{equation}
if simultaneously $\sigma>0$, $h(0)=0$ and $h$ is absolutely continuous,
and $I_\sigma(h)=\infty$ otherwise.
\end{theorem}
The following corollary gives simplified conditions for the MDP principle,
which will be verified in several examples later on.
\begin{corollary}\label{cor1}
Assume that $\|X_0\|_\infty< \infty$ and that $X_0$ is ${\mathcal
F}_0$--measurable. In addition, assume that
\begin{equation}\label{bis}
\sum_{n=1}^{\infty }n^{-1/2}\Vert
\mathbf{E}(X_{n}|\mathcal{F}_{0})\Vert _{\infty }<\infty
\,,
\end{equation}
and that for all $i,j\geq 1$,
\begin{equation}
\lim_{n\rightarrow \infty }\Vert \mathbf{E}(X_{i}X_{j}|\mathcal{F}_{-n})
-\mathbf{E}(X_{i}X_{j})\Vert _{\infty }=0\, .  \label{mix}
\end{equation} Then the conclusion of Theorem \ref{mainresult}
holds with $\sigma^2=\sum_{k \in {\mathbf Z} } {\mathbf E}
(X_0X_k)$.
\end{corollary}
$\quad \;$ The next theorem allows to deal with non-adapted
sequences and it provides additional applications. Let
$\mathcal{F}_{-\infty} = \bigcap_{n \geq0} \mathcal{F}_{-n}$ and
${\mathcal {F}}_{\infty} = \bigvee_{k \in {\mathbf Z}} {\mathcal
{F}}_{k}$.
\begin{theorem}\label{projresult}
Assume that  $\|X_0\|_\infty<\infty$, $\E(X_0|\F_{-\infty})=0$
almost surely, and $X_0$ is $\F_{\infty}$-measurable. Define the
projection operators by
$P_{j}(X)=\mathbf{E}(X|\mathcal{F}_{j})-\mathbf{E}(X|\mathcal{F}_{j-1})\;.$
Suppose that (\ref{mix}) holds and that
\begin{equation}
\sum_{j\in \mathbf{Z}}\Vert P_{0}(X_{j})\Vert _{\infty }<\infty \,.
\label{projcond}
\end{equation}
Then the conclusion of Theorem \ref{mainresult}
holds with $\sigma^2=\sum_{k \in {\mathbf Z} } {\mathbf E}
(X_0X_k)$.
\end{theorem}

\subsection{Simple examples, comments and extensions}

\begin{comment} \begin{rm}
{\bf $\phi$-mixing sequences}. Recall that if $Y$ is a random
variable with values in a Polish space ${\mathcal Y}$ and if
${\mathcal M}$ is a $\sigma$-field, the $\phi$-mixing coefficient
between ${\mathcal M}$ and $\sigma(Y)$ is defined by
\begin{equation}\label{defphi}
\phi({\mathcal M}, \sigma (Y)) =\sup_{A \in {\mathcal B}({\mathcal
Y})} \| {\p}_{Y|{\mathcal M}}(A)-{\p}_Y(A)\|_\infty \, .
\end{equation}
For the sequence $(X_i)_{i \in {\Z}}$ and positive integer $m$ , let
$\phi_m(n)= \sup_{i_m>\ldots>i_1 \geq n} \phi({\mathcal M}_0, \sigma
(X_{i_1},\ldots,X_{i_m}))$ and let
$\phi(k)=\phi_\infty(k)=\lim_{m\to\infty}\phi_m(k)$ be the usual
$\phi$--mixing coefficient. It follows from Corollary \ref{cor1}
that if  the variables are bounded, the conclusion of Theorem
\ref{mainresult} holds as soon as
\begin{equation} \label{condphi2} \sum_{k
>0} k^{-1/2}\phi_1(k)< \infty \;\; \mbox{and} \;\;
\lim_{k\to\infty}\phi_2(k)=0 \, .
\end{equation} The condition (\ref{condphi2})
improves on the one imposed by Gao (1996), that is $\sum_{k >
0}\phi(k)< \infty$, to get the MDP for bounded random variables (see
his  Theorem 1.2).
\end{rm}
\end{comment}

\begin{comment} \label{comLIL} \begin{rm}
{\bf Application to the functional LIL}. Since the variables are
bounded, under the assumptions of Theorem \ref{mainresult} or of
Theorem \ref{projresult}, the MDP also holds in $C[0,1]$ for the
Donsker process $$D_{n}(t)=W_n(t) + n^{-1/2}(nt -[nt])X_{[nt]+1} \,
.$$Hence, if $\sigma^2 > 0$, it follows from the proof of Theorem
1.4.1 in Deuschel and Stroock (1989), that the process
\begin{equation} \label{defdl}\{ ( 2 \sigma^2 \log \log n
)^{-1/2}D_n(t) : t \in [0,1] \} \end{equation} satisfies the
functional law of the iterated logarithm. To be more precise, if
${\mathcal S}$ denotes the subset of $C[0,1]$ consisting of all
absolutely continuous functions with respect to the Lebesgue measure
such that $h(0)=0$ and $\int_0^1 (h'(t))^2 dt \leq 1$, then the
process defined in (\ref{defdl}) is relatively compact with a.s.
limit set ${\mathcal S}$. In the case of bounded random variables,
we then get new criteria to derive the functional LIL. In
particular, the functional LIL holds for $\phi$-mixing bounded
random variables satisfying (\ref{condphi2}).
\end{rm}
\end{comment}

\begin{comment} \begin{rm}  {\bf Linear processes}. Let $(c_i)_{i \in
{\bf Z}}$ be a sequence of real numbers in $\ell^1(\Z)$ (absolutely
summable). Define  $X_k = \sum_{i \in {\Z}}c_i \varepsilon_{k-i}$
where $(\varepsilon_k)_{k \in {\bf Z}}$ is a strictly stationary
sequence satisfying (\ref{mix}) and (\ref{projcond}). Then, so does
the sequence $(X_k)_{k \in{\bf Z}}$, and the conclusion of Theorem
\ref{projresult} holds. In particular, the result applies if
$\varepsilon_0$ is ${\mathcal F}_0$-measurable,  $\E(\varepsilon_1 |
{\mathcal F}_0)=0$ and
\begin{equation*}
 \lim_{n \rightarrow \infty}\Vert {\mathbf E}
(\varepsilon_0^2 |  {\mathcal F}_{-n})  -{\mathbf E}
(\varepsilon_0^2) \Vert_{\infty} = 0 \ .
\end{equation*}
\end{rm}
\end{comment}

\begin{comment} \begin{rm} {\bf Non--mixing in the ergodic sense example.}
The following simple
example shows that Theorem \ref{mainresult} is applicable to
non-mixing in the ergodic theoretical sense sequences. Moreover it
covers a strictly larger class of examples than its Corollary
\ref{cor1}. For all ${ k \in {\mathbf Z}}$, let $Q_{k+1} =-Q_k$
where ${\mathbf P}(Q_0=\pm1)=1/2$ and $X_k=Q_k+Y_k$ where $(Y_k)_{k
\in {{\mathbf Z}}}$ is an iid sequence of zero mean and bounded
random variables, independent of $Q_0$. We can easily check that all
the conditions of Theorem \ref{mainresult} hold while the conditions
of Corollary \ref{cor1} are not satisfied.
\end{rm}
\end{comment}

\begin{comment} \begin{rm} {\bf Stationary ergodic martingales that
does not satisfy MDP.} Let $Y_k$ be the stationary discrete Markov
chain with the state space ${\mathbf N}$ and the transition kernel
given by $\mathbf{P}(Y_1=j-1|Y_0=j)=1$ for all $j\geq 1$ and
$\mathbf{P}(Y_1=j|Y_0=0)=\mathbf{P}(\tau=j)$ for $j\in {\mathbf N}$
with $\E(\tau)<\infty$ and $\mathbf{P}(\tau=1)>0$ which implies that
$(Y_k)$ is ergodic. Let $X_k=\xi_k I_{(Y_k\neq 0)}$ where $(\xi_k)$
is an iid sequence independent of $(Y_k)$ and such that
$\mathbf{P}(\xi_k=\pm 1)=1/2$. Then $X_k$ is a stationary ergodic
martingale difference which is also a bounded function of an ergodic
Markov chain. Straightforward computations show that if $\tau$ does
not have a finite exponential moment then there exists a positive
sequence $a_n\to 0$ with $na_n\to \infty$ for which (\ref{mdpkey})
does not hold. Thus the MDP principle is not true in general for the
stationary sequences satisfying (\ref{bis}) without a certain form
of condition (\ref{s2inf}). A similar example was suggested in
Djellout (2002, Remark 2.6).
\end{rm}
\end{comment}

\begin{comment} \begin{rm} {\bf On $\var(S_n)$ and Theorem
\ref{mainresult}.} Note that if $\sum_{n=1}^{\infty} {n^{-3/2}}\Vert
\E (S_n | {\cal F}_{0} ) \Vert_2  < \infty$, then, by Peligrad and
Utev (2005)
\begin{equation*}
\lim_{n\to\infty} \frac{\var(S_n)}{n} = \sigma ^{2}=
\mathbf{E}(X_{1}^{2})+\sum_{j=0}^{\infty }2^{-j}\mathbf{E}
(S_{2^{j}}(S_{2^{j+1}}-S_{2^{j}})) \, .
\end{equation*}
On the other hand, we shall prove later on that Condition (\ref{mw})
along with (\ref{mix}) are sufficient for the validity of
(\ref{s2inf}). Therefore the conclusion of Theorem \ref{mainresult}
holds under (\ref{mw}) and (\ref{mix}) with $\sigma ^{2}$ identified
in this remark.
\end{rm}
\end{comment}

\begin{comment} \begin{rm} {\bf Sequences that are not strictly stationary}. The proof of
Theorem \ref{projresult} is based on the exponential inequality from
Lemma \ref{exp2}, that was established without stationarity
assumption. Therefore, Theorem \ref{projresult} admits various
extensions to non-stationary sequences. The following slight
generalization is motivated by the fixed design regression problem
 $Z_k=\theta
q_k +X_k$, where the fixed design points are of the form $q_k=1/g(k/n)$,
the error process $X_k$ is a stationary sequence and we analyze the
error of the estimator $\hat{\theta}=n^{-1}\sum_{k=1}^n Z_k g(k/n)$.
If $\{X_{i}\}_{i \in {\bf Z}}$ satisfies the conditions of Theorem
\ref{projresult}, and if $g$ is a Lipschitz function, then the
process $W_n=\{n^{-1/2}\sum_{i=1}^{[nt]} g(i/n) X_i, t \in [0,1]\}$
satisfies (\ref{mdpkey}) with the good rate function $J(\cdot)$
defined by
\begin{equation*}
J(h) = \frac{1}{2\sigma^2}\int_{0}^{1}\Big (
\frac{h^{\prime}(u)}{g(u)} \Big )^2 du \, ,
\quad\mbox{where}\quad \sigma^2=\sum_{k \in {\mathbf Z} } {\mathbf E} (X_0X_k).
\end{equation*}
The proof of this result is omitted. It can be done by following
the proof of Theorem \ref{projresult}. To be more precise, we
start by proving the MDP for the process $ \overline W_n(t) =
n^{-1/2}\sum_{i=1}^{[nv^{-1}(t)]} g(i/n) X_i$ where $v(t)  =
\sigma^2 \int_{0}^{t}g^{2}(x)dx$. For $\overline W_n( . ) $, the
rate function is $I_{\sigma}(.)$ as in Theorem \ref{mainresult}.
To go back to the process $W_{n}(.)$, use the change-of-time
$W_{n}=\overline W_n \circ v$.

\end{rm}
\end{comment}

\section{Applications}
In this Section we present applications to functions of
$\phi$-mixing processes, contracting Markov chains, expanding maps
of the interval and symmetric random walks on the circle. The proofs
are given in Section 4.

\subsection{Functions of $\phi$-mixing sequences}\label{funclin}

In this section, we are partly motivated by Djellout et al. (2006,
Theorem 2.7), who have proved   the MDP  for
\begin{equation}\label{Djetal}X_k=f(Y_k, \dots, Y_{k - \ell}) - {\mathbf E}(f(Y_k,
\dots, Y_{k - \ell}))\;\mbox{ where }\; Y_k= \sum_{i \in {\mathbf
Z}} c_i
 \varepsilon_{k-i}
 \end{equation}
In their Theorem 2.7, Djellout et al. (2006) assume that\\
(i) $(\varepsilon_i)_{i \in {\mathbf Z}}$ is an iid sequence;\\
(ii) (condition on $c_i$) the spectral density of $Y_k$ is continuous on $[-\pi, \pi[$;\\
(iii) (condition on $\varepsilon_0$) $\varepsilon_0$ satisfies the
so-called LSI condition, which implies  that $\E(\exp (\delta
\varepsilon_0^2))<\infty $ for some positive $\delta$, and that the
distribution $\varepsilon_0$ has an absolutely continuous component
 with respect to the Lebesgue measure with a
 strictly positive density on the  support of $\mu$
(see their condition (2.1));\\
(iv) (condition on  $f$) the functions $\partial_{i} f $ are
Lipschitz for $i=0, \dots, \ell$.

By applying our main results, we derive the Propositions
\ref{proplip1} and \ref{proplip2} stated below. In the case where
$X_k$ is given by (\ref{Djetal}), the  Proposition \ref{proplip1}
will allow us to obtain the MDP for a  large class of functions.
However, we require a stronger condition than (ii), that is we
assume that the sequence $(c_i)_{i \in \mathbf Z}$ is in
$\ell_1(\Z)$, and instead of (iii), we suppose that $\varepsilon_0$
takes its values in  some compact intervall $[a,b]$ (this assumption
cannot be compared  to the LSI condition (iii)). Our method allows
to link the regularity of $f$ to the behavior of the coefficients
$(c_i)_{i \in \mathbf Z}$ (in that case, the condition
(\ref{modulus}) given below means that $\sum_{i \in \Z} w_j
(2(b-a)|c_i|)< \infty$ for any $j=0, \ldots, \ell$, where $w_j$ is
the modulus of continuity of $f$ with respect to the $j$-th
coordinate). In addition, our innovations maybe dependent: more
precisely, $(\varepsilon_i)_{i \in {\mathbf Z}}$ is assumed to be a
stationary $\phi$--mixing sequence. \vskip3pt

We now describe our general results. Let $(\varepsilon_i)_{i \in
{\mathbf Z}}=(\varepsilon_0 \circ T^i)_{i \in {\mathbf Z}}$ be a
stationary sequence of $\phi$-mixing random variables with values in
a subset $A$ of a Polish space ${\mathcal X}$. Starting from the
definition (\ref{defphi}), we denote by $\phi_\varepsilon (n)$ the
coefficient $\phi_\varepsilon (n)= \phi (\sigma(\varepsilon_i, i
\leq 0), \sigma( \varepsilon_i, i\geq n))$.

Our first result is for non-adapted sequences, that is satisfying
the  representation (\ref{HX}) below. Let $H$ be a function from
$A^{\mathbf Z}$ to ${\mathbb R}$ satisfying the condition
$$
C(A): \quad \text{for any $x, y$ in $A^{\mathbb Z}$},\quad
|H(x)-H(y)| \leq \sum_{i\in {\mathbf Z}} \Delta_i {\bf 1}_{x_i\neq
y_i}, \quad \text{where $\displaystyle \sum_{i \in {\mathbf
Z}}\Delta_i < \infty$,}
$$
 Define the stationary sequence
$X_k=X_0\circ T^k$ by
\begin{equation}\label{HX}
X_k=H((\varepsilon_{k-i})_{i \in {\mathbf
Z}})-{\E}(H((\varepsilon_{k-i})_{i \in {\mathbf Z}}))\, .
\end{equation}
Note that $X_k$ is bounded  in view of $C(A)$.
\begin{proposition}\label{proplip1}
Let $(X_k)_{k \in {\mathbb Z}}$ be defined by (\ref{HX}), for a
function $H$ satisfying $C(A)$. If $\sum_{k>0} \phi_\varepsilon(k)$
is finite, then the conclusion of Theorem \ref{mainresult} holds
with $\sigma^2=\sum_{k \in {\mathbf Z} } {\mathbf E} (X_0X_k)$.
\end{proposition}

For adapted sequences, that is satisfying the representation
(\ref{HXC}) below, we can assume that $H$ satisfies another type of
condition. Let $H$ be a function from $A^{\mathbf N}$ to ${\mathbb
R}$ satisfying the condition
$$
C'(A): \quad \text{for any $i \geq 0$}, \ \sup_{x \in A^{\mathbf N},
y \in A^{\mathbf N} } |H(x)-H(x^{(i)}y)| \leq R_i,  \quad
\text{where $R_i$ decreases to  0,}
$$
the sequence $x^{(i)}y$ being defined by  $(x^{(i)}y)_j=x_j$ for $j<
i$ and $(x^{(i)}y)_j= y_j$ for $j \geq i$. Define the stationary
sequence $X_k=X_0\circ T^k$ by
\begin{equation}\label{HXC}
X_k=H((\varepsilon_{k-i})_{i \in {\mathbf
N}})-{\E}(H((\varepsilon_{k-i})_{i \in {\mathbf N}}))\, .
\end{equation}
\begin{proposition}\label{proplip2}
Let $(X_k)_{k \in {\mathbb Z}}$ be defined by (\ref{HXC}), for a
function $H$ satisfying $C'(A)$. If
\begin{equation}\label{mixrphi}
\sum_{\ell=1}^\infty R_\ell \sum_{k\geq \ell}
\frac{\phi_\varepsilon(k-\ell)}{\sqrt k}< \infty\, ,
\end{equation}
 then the conclusion of Theorem \ref{mainresult} holds with
$\sigma^2=\sum_{k \in {\mathbf Z} } {\mathbf E} (X_0X_k)$. In
particular, the condition (\ref{mixrphi}) holds as soon as
\begin{enumerate}
\item $\sum_{k>0} \phi_\varepsilon(k)< \infty$ and $\sum_{k>0}
k^{-1/2} R_k < \infty$.
\item $\sum_{k>0} R_k< \infty$ and $\sum_{k>0}
k^{-1/2} \phi_\varepsilon(k) < \infty$.
\end{enumerate}
\end{proposition}
{\bf Application to functions of linear processes.} Assume that
$\varepsilon_i$ takes its values in a compact interval $A=[a, b]$ of
${\mathbb R}$, and let $(c_i)_{i \in {\mathbf Z}}$ be a sequence of
real numbers in $\ell^1(\Z)$. Let $m=\inf_{x \in A^{\mathbf Z}}
\sum_{i \in \Z} c_i x_i$ and $M=\sup_{x \in A^{\mathbf Z}} \sum_{i
\in \Z} c_i x_i$.
 For a function  $f$ from $[m, M]^\Z$ to ${\mathbb R}$, let $w_i$ be the
 modulus of continuity of $f$ with respect to the $i$-th coordinate, that
 is
 $$
  w_i(h)=\sup_{x \in [m, M]^\Z, t \in [m, M], |x_i-t|\leq h}
  |f(x)-f(x^{(i,t)})|\, ,
 $$
 the sequence $x^{(i,t)}$ being defined by  $x^{(i,t)}_j=x_j$ for
$j\neq i$ and $x_i^{(i,t)}= t$.
 Assume that
 $$
\text{for any $x, y$ in $[m,M]^\Z$} \quad |f(x)-f(y)| \leq \sum_{i
\in {\mathbb Z}} w_i(|x_i-y_i|) < \infty \, .
 $$
 Define the random
variables $Y_k= \sum_{i \in {\mathbf Z}} c_i
 \varepsilon_{k-i}$, and let
 \begin{equation}\label{gener}
   X_k= f((Y_{k-i})_{i \in {\mathbb Z}})-\E(f((Y_{k-i})_{i \in {\mathbb
   Z}})\,
 \end{equation}
 (note that (\ref{gener}) is a generalization of (\ref{Djetal})).
 Clearly, $X_k$ may be written as in (\ref{HX}), for a function $H$ from $A^\Z$ to
 ${\R}$. Moreover, $H$ satisfies $C(A)$ with $ \Delta_i \leq
 \sum_{\ell \in \Z} w_\ell(2(b-a)|c_{i-\ell}|)$ provided that
 \begin{equation}\label{modulus}
\sum_{i \in \Z}\sum_{\ell \in \Z} w_\ell(2(b-a)|c_{i}|)< \infty \, .
 \end{equation}
\vskip10pt
 From Proposition \ref{proplip1}, if $\sum_{k>0} \phi_\varepsilon(k)
< \infty$ and if (\ref{modulus}) holds, then the conclusion of
Theorem \ref{mainresult} holds.
 In particular, the condition (\ref{modulus}) holds as soon as there exist
 $(b_i)_{i \in \Z}$ in $\ell^1(\Z)$ and $\alpha$ in $]0, 1]$ such
 that $w_\ell(h)\leq
 b_\ell|h|^\alpha$ and
 $\sum_{i \in \Z} |c_i|^\alpha<\infty$. Two simple examples of such functions
 are:
 \begin{enumerate}
 \item $f(x)= \sum_{i\in \Z} g_i(x_i)$ for some  $g_i$ such that
 $|g_i(x)-g_i(y)|\leq b_i|x-y|^\alpha$ for any $x, y$ in $[m, M]$.
\item $ f(x)=\Pi_{i=p}^q h_i(x_i)$ for some  $h_i$ such
that $|h_i(x)-h_i(y)|\leq K_i|x-y|^\alpha$ for any $x, y$ in $[m,
M]$.
\end{enumerate}

Now, assume that $c_i=0$ for $i<0$, so that
$Y_k= \sum_{i \geq 0} c_i
 \varepsilon_{k-i}$. If $f$ is in fact a
function of $x$ through $x_0$ only, we simply denote by $w=w_0$ its
modulus of continuity over $[m, M]$. In that case
$X_k=f(Y_k)-{\mathbf E}(Y_k)$ may be written as in (\ref{HXC}) for a
function $H$ satisfying $C'(A)$ with  $R_i \leq w(2|b-a|\sum_{k\geq
i}|c_k|)$. From item 1 of Proposition \ref{proplip2}, if $\sum_{k>0}
\phi_\varepsilon(k) < \infty$ and if
\begin{equation}\label{Kac}
\sum_{n \geq 1} n^{-1/2}w\Big(2|b-a| \sum_{k \geq n} |c_k|\Big )
 < \infty \, ,
\end{equation}
then the conclusion of Theorem \ref{mainresult} holds. In
particular, if $|c_i|\leq C \rho^i$ for some $C>0$ and $\rho \in ]0,
1[$,  the condition (\ref{Kac}) holds as soon as:
$$
  \int_0^1 \frac{w(t)}{t \sqrt{|\log t|}} dt < \infty \, .
$$
Note that this condition  is satisfied as soon as $w(t)\leq
D|\log(t)|^{-\gamma}$ for some $D>0$ and some $\gamma> 1/2$. In
particular, it is satisfied if $f$ is $\alpha$-H\"older for some
$\alpha \in ]0, 1]$.

\subsection{Contracting Markov chains}
Let $(Y_n)_{n \geq 0}$ be a stationary Markov chain of bounded
random variables  with invariant measure $\mu$ and transition kernel
$K$. Denote by $\|\cdot \|_{\infty, \mu}$ the essential supremum
norm with respect to $\mu$. Let $\Lambda_1$ be the set of
$1$-Lipschitz functions. Assume that the chain satisfies the two
following  conditions:
\begin{eqnarray}
\text{there exist $C>0$ and $\rho \in ]0, 1[$ such that} \quad
\sup_{g \in {\Lambda_1}} \Vert K^n (g)
- \mu (g) \Vert_{\infty, \mu}\leq C \rho^n \, ,\label{contrac1}\\
\text{for any $f,g \in \Lambda_1$ and any $m\geq 0$} \quad \lim_{n
\rightarrow \infty} \|K^n(fK^m(g))-\mu(f K^m(g))\|_{\infty, \mu}=0
\, . \label{contrac2}
\end{eqnarray}
We shall see in the next proposition that if (\ref{contrac1}) and
(\ref{contrac2}) are satisfied, then the MDP holds in $D[0, 1]$ for
the sequence
\begin{equation}\label{defxn}
X_n=f(Y_n)-\mu(f)
\end{equation}
 as soon as  the function $f$ belongs to the class ${\cal L}$ defined below.
\begin{definition}
Let ${\cal L}$ be the class of functions $f$ from ${\bf R}$ to ${\bf
R}$ such that $ | f(x) - f(y) | \leq c( |x-y|) $, for some concave
and non decreasing function $c$ satisfying
\begin{equation}\label{intc}
\int_0^1 \frac{c(t)}{t \sqrt{|\log t |}} dt < \infty \, .
\end{equation}
\end{definition}
Note that (\ref{intc}) holds if $c(t)\leq D|\log(t)|^{-\gamma}$ for
some $D>0$ and some $\gamma> 1/2$. In particular, ${\mathcal L}$
contains   the class of functions from $[0,1]$ to ${\bf R}$ which
are $\alpha$-H\"older for some $\alpha \in ]0, 1]$.
\begin{proposition}\label{contracting}
Assume that the stationary Markov chain $(Y_n)_{n \geq 0}$ satisfies
(\ref{contrac1}) and (\ref{contrac2}),  and let $X_n$ be defined by
(\ref{defxn}). If $f$ belongs to ${\cal L} $, then the conclusion of
Theorem \ref{mainresult} holds with
$$
\sigma^2=\sigma^2(f)= \mu((f -\mu(f))^2) + 2 \sum_{n > 0} \mu(K^n(f)
\cdot (f-\mu(f))) \, .
$$
\end{proposition}
The proof of this proposition is based on the following lemma which
has interest in itself.
\begin{lemma} \label{rhoarith}
Let  $u_n = \sup_{g \in {\Lambda_1}} \Vert K^n (g) - \mu (g)
\Vert_{\infty, \mu}$. Let $f$ be a function from $\bf R$ to ${\bf
R}$ such that $ |f(x) - f(y)| \leq c(|x-y|) $ for some concave and
non decreasing function $c$. Then
$$
\|K^n(f)-\mu(f)\|_{\infty, \mu} \leq c(u_n) \, .
$$
\end{lemma}
\begin{remark}
If $u_n \leq C \rho^n$ for a $C > 0$ and $\rho \in ]0,1[$, and if
$c(t) \leq D| \log (t) |^{- \gamma}$ for $D >0$ and $\gamma >0$,
then
$$
\|K^n(f)-\mu(f)\|_{\infty, \mu} = O (n^{- \gamma}) \, .
$$
\end{remark}

We now give two conditions under which (\ref{contrac1}) and
(\ref{contrac2}) hold. Let $[a, b]$ be a compact interval in which
lies the support of $\mu$. For a Lipschitz function $f$, let
$\text{Lip}(f)=\sup_{x, y \in [a, b]}|f(x)-f(y)|/|x-y|$. The chain
is said to be Lipschitz contracting if there exist $\kappa>0$ and
$\rho \in ]0, 1[$ such that
\begin{equation}\label{lip}
  \text{Lip}(K^n(f)) \leq \kappa \rho^n \text{Lip}(f)\, .
\end{equation}
Let $BV$ be the class of bounded variation functions from $[a,b]$ to
${\R}$. For any $f \in BV$, denote by $\|df \|$ the total  variation
norm of the measure $df$: $\|df\|=\sup\{\int g df, \|g\|_\infty \leq
1\}$. The chain is said to be to be $BV$-contracting if there exist
$\kappa
>0$ and $\rho \in [0, 1[$ such that
\begin{equation}\label{expan}
     \|d K^n(f)\| \leq \kappa \rho^n \|df \| \, .
\end{equation}
It is easy to see that if either (\ref{lip}) or (\ref{expan}) holds,
then (\ref{contrac1}) and (\ref{contrac2}) are satisfied (to see
that the condition (\ref{expan}) implies (\ref{contrac2}), it
suffices to note that it implies the same property for two $BV$
functions $f, g$ (see (\ref{2points})), and that any Lipshitz
function from $[a, b]$ to $\R$ can be uniformly approximated by $BV$
functions).

\medskip

\noindent{\bf Application to iterated random functions.} The
stationary bounded Markov chain $(Y_n)_{n \geq 0}$ with transition
kernel $K$  is one-step Lipschitz contracting if there exists $\rho
\in ]0, 1[$ such that
$$
 \text{Lip}(K(f))\leq  \rho \text{Lip}(f) \, .
$$
Note that if $K$ is one-step Lipschitz contracting then (\ref{lip})
obviously holds with $\kappa=1$. The one-step contraction is a very
restrictive assumption. However, it is satisfied if $Y_n=F(Y_{n-1},
\varepsilon_n)$ for some iid sequence $(\varepsilon_i)_{i >0}$
independent of $Y_0$, and some function $F$ such that
\begin{equation}\label{1step}
\|F(x, \varepsilon_1)-F(y, \varepsilon_1)\|_1 \leq \rho |x-y| \quad
\text{for any $x,y$ in $\R$.}
\end{equation}
\begin{remark}
Under a more restrictive condition on $F$ than (\ref{1step}), namely
\begin{equation}\label{Flip}
|F(x,z)-F(y,t)|\leq \rho|x-y| + L|z-t| \, ,
\end{equation}
Delyon et al (2006) have proved the MDP for $X_n=f(Y_n)-\mu(f)$ when
$f$ is a lipschitz function. In their paper, the chain is not
assumed to be bounded. It is only assumed that $\E(e^{\delta
\varepsilon_1})< \infty $ for some $\delta>0$, which implies the
same property for $X_1$ (for a smaller $\delta$) by using the
inequality (\ref{Flip}).
\end{remark}

\noindent{\bf Application to expanding maps.}
Let $T$ be a map from $[0, 1]$ to $[0, 1]$ preserving a probability
$\mu$ on $[0, 1]$,  and let
\begin{eqnarray*}
X_k = f \circ T^{n-k+1}- \mu(f)  \, ,\;  W_{n}(t) = W_n(f,t) =
n^{-1/2}\sum_{i=1}^{[nt]} (f \circ T^{n-i+1}- \mu(f))
\end{eqnarray*}
Define  the Perron-Frobenius operator $K$ from $L^2([0, 1], \mu)$ to
$L^2([0, 1], \mu)$ $via$ the equality
\begin{equation}\label{eq}
   \int_0^1 (Kh)(x) f(x) \mu (dx) = \int_0^1 h(x) (f \circ T)(x)
   \mu(dx)
   \, .
\end{equation}
The map $T$ is said to be $BV$-contracting if its Perron-Frobenius
operator is $BV$-contracting, that is satisfies (\ref{expan}). As a
consequence of Proposition \ref{contracting}, the following
corollary holds.
\begin{corollary}\label{proexp}
If $T$ is $BV$-contracting, and if $f$ belongs to $BV \cup {\cal L}
$, then the conclusion of Theorem \ref{mainresult} holds with
$$
\sigma^2=\sigma^2(f)= \mu((f -\mu(f))^2) + 2 \sum_{n > 0} \mu(f
\circ T^n \cdot (f-\mu(f))) \, .
$$
\end{corollary}

Let us present a large class of $BV$-contracting maps. We shall say
that  $T$ is uniformly expanding if it  belongs to the class
${\mathcal C}$ defined in Broise (1996), Section 2.1 page 11. Recall
that if $T$ is uniformly expanding, then there exists a probability
measure $\mu$ on $[0, 1]$, whose density $f_\mu$ with respect to the
Lebesgue measure is a bounded variation function, and such that
$\mu$ is invariant by $T$. Consider now the more restrictive
conditions:
\begin{enumerate}
\item[(a)] $T$ is uniformly expanding.
\item[(b)] The invariant measure $\mu$ is unique and $(T, \mu)$ is
mixing in the ergodic-theoretic sense.
\item[(c)] $\displaystyle \frac{1}{f_\mu}{\bf 1}_{f_\mu>0}$ is a
bounded variation function.
\end{enumerate}
Starting from Proposition 4.11 in Broise (1996), one can prove that
if $T$ satisfies the assumptions (a), (b) and (c) above, then it is
$BV$ contracting (see for instance Dedecker and Prieur (2007),
Section 6.3). Some well known examples of maps satisfying the
conditions (a), (b) and (c) are:
\begin{enumerate}
\item $T(x)= \beta x -[\beta x]$ for $\beta> 1$. These maps are called
$\beta$-transformations.
\item  $I$ is the finite union of disjoint intervals $(I_k)_{1 \leq k
\leq n}$, and $T(x)=a_kx +b_k$ on $I_k$, with $|a_k|>1$.
\item $T(x)=a(x^{-1}-1)-[a(x^{-1}-1)]$ for some $a>0$. For $a=1$, this
transformation is known as the Gauss map.
\end{enumerate}

\begin{remark}
The case where $f(x)=x$ (that is $X_n=T^n- \mu(T)$)  has already
been considered by Dembo and Zeitouni (1997). However, in this
paper, the assumptions on $T$ are more restrictive than the
assumptions (a), (b) and (c) above. In particular, they assume that
there is a finite partition $(I_j)_{1\leq j \leq m}$ of $[0, 1]$ on
which $T$ restricted to $I_k$ is $C^1$ and $\inf_{x \in I_k}
|T'(x)|>1$, so that their result does not cover the case of the
Gauss map (Example 3 above).
\end{remark}

\subsection{Symmetric random walk on the circle}
$\quad \;$ Let $K$ be the Markov kernel defined by $$Kf (x) =
\frac{1}{2}( f (x+a) + f (x-a) )$$ on the torus $\R/\Z$, with $a$
irrational in $[0,1]$. The Lebesgue-Haar  measure $m$  is the unique
probability which is invariant by $K$. Let $(\xi_i)_{i\in \Z}$ be
the stationary Markov chain with transition kernel $K$  and
invariant distribution $m$. Let
\begin{equation}\label{defSnf}
X_k=f(\xi_k)-m(f)\;,\; W_{n}(t) = W_n(f,t) =
n^{-1/2}\sum_{i=1}^{[nt]}(f(\xi_i)-m(f)) \, .
\end{equation}
From Derriennic and Lin (2001), Section 2, we know that the central
limit theorem holds for $n^{-1/2} W_n (f, 1)$ as soon as the series
of covariances
\begin{equation}\label{sericov}
\sigma^2(f) = m((f-m(f))^2) + 2 \sum_{n>0} m( f K^n(f-m(f)))
\end{equation}
is convergent, and that the limiting distribution is ${\mathcal
N}(0,\sigma^2(f))$. In fact the convergence of the series in
(\ref{sericov}) is equivalent to
\begin{equation}\label{Paroux}
 \sum_{k\in{\Z}^*}   \frac{|\hat f (k)|^2}{d (ka, {\Z})^2} < \infty
 \, ,
\end{equation}
where $\hat f(k)$ are the Fourier coefficients of $f$. Hence, for
any irrational number $a$, the criterion (\ref{Paroux}) gives a
class of function $f$ satisfying the central limit theorem, which
depends on the sequence $((d (ka, {\Z}))_{k \in{\Z}^*}$. Note that a
function $f$ such that
\begin{equation}\label{debilmental}
\liminf_{k \rightarrow \infty} k |\hat f(k)|>0 \, ,
\end{equation}
does not satisfies (\ref{Paroux}) for any irrational number $a$.
Indeed, it is well known from the theory of continued fraction that
if $p_n/q_n$ is the $n$-th convergent of $a$, then $|p_n-q_n
a|<q_n^{-1}$, so that $d(ka, \Z)<k^{-1}$ for an infinite number of
positive integers $k$. Hence, if (\ref{debilmental}) holds, then
$|\hat f(k)|/d(ka, \Z)$ does not even tend to zero as $k$ tends to
infinity.

Our aim in this section is to give conditions on $f$ and on the
properties of the irrational number $a$ ensuring that the MDP holds
in $D[0, 1]$.
\begin{eqnarray}\label{badly} && \mbox{$a$ is said to be badly approximable
by rationals if for any positive $\varepsilon$,}\\
&&\mbox{the inequality $d(ka, \Z) < |k|^{-1-\varepsilon}$ has only finitely many
solutions for $k \in \Z$.}\nonumber
\end{eqnarray}

\noindent From Roth's theorem the algebraic numbers are badly
approximable (cf. Schmidt (1980)). Note also that the set of badly
approximable numbers in $[0,1]$ has Lebesgue measure $1$.

In Section 5.3 of Dedecker and Rio (2006), it is proved that the
condition (\ref{Paroux}) (and hence the central limit theorem for
$n^{-1/2} W_n (f, 1)$) holds for any badly approximable number $a$
as soon as
\begin{equation}\label{antideb}
\sup_{k\not= 0} |k|^{1+\varepsilon} |\hat f (k)| < \infty \quad
\text{for some positive $\varepsilon$.}
\end{equation}
Note that, in view of (\ref{debilmental}), one cannot take
$\varepsilon=0$ in  the condition (\ref{antideb}).

 In fact, for badly approximable numbers, the condition  (\ref{antideb}) implies  also the
 MDP in $D[0,1]$:
\begin{proposition}\label{circle}
Suppose that $a$ is badly approximable by rationals, i.e satisfies
(\ref{badly}). If the function $f$ satisfies (\ref{antideb}), then
the conclusion of Theorem \ref{mainresult} holds with
$\sigma^2=\sigma^2(f)$.
\end{proposition}
Note that, under the same conditions, the process $\{W_n(f,t), t \in
[0, 1] \}$ satisfies the weak invariance principle in $D[0, 1]$.
Indeed, to prove Proposition \ref{circle}, we show that the
conditions of Corollary \ref{cor1} are satisfied, but these
conditions imply the weak invariance principle (see for instance
Peligrad and Utev (2005)). From Comment \ref{comLIL}, we also infer
that the Donsker process defined in (\ref{defdl}) satisfies the
functional law of the iterated logarithm.

\section{Proofs}

Since the proofs of our results are mainly based on some exponential
bounds for the deviation probability of the maximum of the partial
sums for dependent variables, we present these inequalities, which
have interest in themselves.

\subsection{Exponential bounds for dependent variables}
$\quad \;$
We state first the exponential bound from Proposition 2 in Peligrad,
Utev and Wu (2007) that we are going to use in the proof of the main theorem.

\begin{lemma}\label{puw}
Let $(X_i)_{ i \in {\mathbb Z}}$ be a stationary sequence of random
variables adapted to the filtration $(\mathcal{F}_i)_{ i \in {\mathbb Z}}$.
Then
\begin{equation*}
{\mathbf P}(\max_{1\leq i\leq n}|S_{i}|\geq t)\leq4\sqrt{e}\exp \big
(-t^{2}/2n \big [\Vert
X_{1}\Vert_{\infty}+80\sum_{j=1}^{n}j^{-3/2}\Vert\mathbf{E}(S_{j}|
\mathcal{F}_{0})\Vert_{\infty} \big ]^{2} \big )
\end{equation*}
\end{lemma}

In the next lemma, we bound the maximal exponential moment of the
stationary sequence by using the projective criteria.

\begin{lemma}
\label{exp2} Let $\{Y_{k}\}_{k\in\Z }$ be a sequence of random
variables such that for all $j$, $\E(Y_j|\F_{-\infty})=0$ almost
surely and $Y_j$ is $\F_{\infty}$-measurable. Define the
projection operators by
$P_{j}(X)=\mathbf{E}(X|\mathcal{F}_{j})-\mathbf{E}(X|\mathcal{F}_{j-1})\;.$
Assume that
\begin{equation*}
\Vert P_{k-j}(Y_k)\Vert_\infty\leq p_j\quad\mbox{and}\quad
D:=\sum_{j=- \infty}^{\infty }p_j<\infty
\end{equation*}
Let $\{g_{k}\;,\;k\in\N\}$ be a sequence of numbers and define,
\begin{equation*}
S_{k}=\sum_{i=1}^{k}g_i Y_{i}\;,\;M_{k}=\max_{1\leq j\leq k}S_{j} \ ,\
G_n^2=\sum_{i=1}^n g_i^2
\end{equation*}
Then,
\begin{equation*}
\mathbf{E}\exp (tM_{n})\leq 4\exp (\tfrac{1}{2}G_n^2 D^{2}t^{2})
\, .
\end{equation*}
In particular,
\begin{equation*}
\mathbf{P}\Big(\max_{1\leq k\leq n}|S_{k}|\geq x)\leq 8\exp
\Big(-\frac{x^{2}}{2G_n^2 D^{2}}\Big) \, .
\end{equation*}
\end{lemma}
{\bf Proof}. Start with the decomposition
\begin{equation*}
Y_{k}=\sum_{j=- \infty}^{\infty }P_{k-j}(Y_{k})=\sum_{j=-
\infty}^{\infty }b_{j}P_{k-j}(Y_{k})/b_{j}
\end{equation*}
where $b_{j}=p_j/D\geq \Vert P_{k-j}(Y_{k})\Vert _{\infty }/D$,
 for any $j \in {\mathbf {Z}}$. Then
\begin{equation*}
S_{m}=\sum_{j=- \infty}^{\infty
}b_{j}\sum_{k=1}^{m}g_k P_{k-j}(Y_{k})/b_{j} \, .
\end{equation*}
Thus,
\begin{equation*}
M_{n}\leq \sum_{j=- \infty}^{\infty }b_{j}\max_{1\leq m\leq
n}\sum_{k=1}^{m}P_{k-j}(g_kY_{k})/b_{j}=:\sum_{j=- \infty}^{\infty
}b_{j}M_{n}^{(j)}
\end{equation*}
where $M_{n}^{(j)}$ denotes $\max_{1\leq m\leq
n}\sum_{k=1}^{m}g_kP_{k-j}(Y_{k})/b_{j}.$

Since $\exp (x)$ is convex and non--decreasing and $b_{j}\geq 0$ with
$\sum_{j \in {\mathbf {Z}}}b_{j}=1$,
\begin{equation*}
\mathbf{E}\exp (tM_{n})\leq \mathbf{E}\exp \Big(\sum_{j=-
\infty}^{\infty }b_{j}tM_{n}^{(j)}\Big)\leq \sum_{j=-
\infty}^{\infty }b_{j}\mathbf{E}\exp (tM_{n}^{(j)}) \, .
\end{equation*}
Consider the martingale difference $U_{k} = g_kP_{k-j}(Y_{k})/b_{j}$,
$j=1, \ldots ,n$. Since the variables $Z_k=\exp(t (U_1+ \cdots
+U_k)/2)$ form a submartingale, Doob's inequality yields
\begin{equation*}
\mathbf{E}\exp (tM_{n}^{(j)})=\mathbf{E}\Big(\max_{1\leq k\leq n}Z_{k}^{2}\Big)\leq
4\mathbf{E}Z_{n}^{2}=4\mathbf{E}\exp (t(U_1+ \cdots +U_n))\, .
\end{equation*}
Applying Azuma's inequality to the right-hand side,  and noting that
\begin{equation*}
\Vert U_{k}\Vert _{\infty }=|g_k|\Vert P_{k-j}(Y_{k})\Vert _{\infty
}/b_{j}\leq |g_k|D \, ,
\end{equation*}
we infer that
\begin{equation*}
\mathbf{E}\exp (tM_{n}^{(j)})\leq 4\exp (\tfrac{1}{2}G_n^2D^{2}t^{2}) \,.
\end{equation*}
Since $\sum_{j\in {\mathbf {Z}}}b_{j}=1$, we obtain that
\begin{equation*}
\mathbf{E}\exp (tM_{n})\leq \sum_{j\in {\mathbf {Z}}}b_{j}4\exp
(\tfrac{1}{2}G_n^2D^{2}t^2)= 4\exp(\tfrac{1}{2}G_n^2D^{2}t^2) \, .
\end{equation*}
Next, to derive the one--sided probability inequality we use the
exponential bound with $t=x/(G_n^2D^{2})$, so
\begin{equation*}
\mathbf{P}(M_{n}\geq x)\leq \mathbf{E}\exp (tM_{n})\exp (-tx)=4\exp
\Big(- \frac{x^{2}}{2G_n^2D^{2}}\ \Big) \, .
\end{equation*}
Finally, to derive the two--sided inequality we observe that the
stationary sequence $\{-Y_{j}\}$ also satisfies the conditions of
the lemma. The proof is complete. \quad $\diamond $ \vskip10pt

The next technical lemma provides an exponential bound for any
random vector plus a correction in terms of conditional
expectations (see also Wu, 1999).

\begin{lemma}
\label{techlemma1} Let $\{X_{i}\}_{1\leq i\leq n}$ be a vector of
real random variables adapted to the filtration
$\{\mathcal{F}_{n}\}_{n\geq 1}$. Denote $B$=$\sup_{1\leq i\leq
n}\|X_{i}\|_{\infty }$. Then, for all $\delta
>0$ and $c$ a natural number with $cB/n\leq $ $\delta /2$, we have
\begin{equation}
\mathbf{P} (\max_{1\leq i\leq
n}|\frac{1}{n}\sum_{u=1}^{i}X_{u}|\geq \delta )\leq 2\exp
(-\frac{\delta ^{2}n}{64B^{2}c})+ \mathbf{P} (\sup_{1\leq i\leq
\lbrack n/c]}|\frac{1}{c}\sum_{j=(i-1)c+1}^{ic} \mathbf{E} (X_{j}|
\mathcal{F}_{(i-1)c})|\geq \frac{\delta }{4})  \label{ineq}
\end{equation}
\end{lemma}
{\bf Proof of Lemma \ref{techlemma1}} \ Let $c$ be a
fixed integer and $k=[n/c]$ (where, as before, $[x]$ denotes the integer part of
$x $). The initial step of the proof is to divide the variables in consecutive
blocks of size $c$ and to average the variables in each block
\begin{equation*}
Y_{i,c}=\frac{1}{c}\sum_{j=(i-1)c+1}^{ic}X_{j}\,,\text{ }i\geq 1\,.
\end{equation*}

Then, for all $1\leq i\leq k$ we construct the martingale,
\begin{equation*}
M_{i,c}=\sum_{j=1}^{i}(Y_{j,c}-\mathbf{E}(Y_{j,c}|
\mathcal{F}_{(j-1)c}))=\sum_{j=1}^{i}D_{j,c}\,
\end{equation*}
and we use the decomposition
\begin{eqnarray*}
\mathbf{P}(\max_{1\leq j\leq n}|\frac{1}{n}\sum_{u=1}^{j}X_{u}|
&\geq &\delta )\leq \mathbf{P}(\max_{1\leq i\leq
k}|\frac{1}{k}\sum_{j=1}^{i}Y_{j,c}|\geq \delta -
\frac{cB}{n})\leq \mathbf{P}(\max_{1\leq i\leq
k}|\frac{1}{k}\sum_{j=1}^{i}Y_{j,c}|\geq
\delta /2) \\
&\leq &\mathbf{P}(\max_{1\leq i\leq k}\frac{1}{k}|M_{i,c}|\geq
\delta / 4 )+\mathbf{P}(\max_{1\leq i\leq
k}\frac{1}{k}|\sum_{j=1}^{i}\mathbf{E}(Y_{j,c}|
\mathcal{F}_{(j-1)c})|\geq \delta / 4) \\
&\leq &\mathbf{P} (\max_{1\leq i\leq k}\frac{1}{k}|M_{i,c}|\geq
\delta / 4) + \mathbf{P} (\max_{1\leq j\leq
k}|\mathbf{E}(Y_{j,c}|\mathcal{F}_{(j-1)c})|\geq \delta / 4)\,.
\end{eqnarray*}

Next, we apply Azuma's inequality to the martingale part and obtain,
\begin{equation*}
\mathbf{P}\left( \max_{1\leq i\leq k}|M_{i,c}|\geq \delta k/4\right)
\leq
2\exp (-\frac{\delta ^{2}k^{2}}{32kB^{2}})\leq 2\exp (-\frac{\delta ^{2}n}
{64cB^{2}})
\end{equation*}
which implies that
\begin{equation*}
\mathbf{P}(\max_{1\leq i\leq n}|\frac{1}{n}\sum_{u=1}^{i}X_{u}|\geq
\delta )\leq 2\exp (-\frac{\delta
^{2}n}{64B^{2}c})+\mathbf{P}(\max_{1\leq i\leq k}|\mathbf{E}
(Y_{i,c}|\mathcal{F}_{(i-1)c})|\geq \delta / 4)
\end{equation*}\
proving the lemma. \quad $\diamond $

\subsection{Some facts about the moderate deviation principle}
$\quad \;$ This paragraph deals with some preparatory material.
The following theorem is a result concerning the MDP for a
triangular array of martingale differences sequences. It follows
from Theorem 3.1 and Lemma 3.1 of Puhalskii (1994), (see also
Djellout (2002), Proposition 1 and Lemma 2).

\begin{lemma}\label{triangmart}
Let $k_n$ be an increasing sequence of integers going to infinity.
Let $\{D_{j,n}\}_{1\leq j\leq k_{n}}$ be a triangular array of
martingale differences adapted to a filtration
${\mathcal{F}_{j,n}}$. Define the normalized partial sums process
$Z_{n}(t)=n^{-1/2}\sum_{i=1}^{[k_{n}t]}D_{i,n}$. Let $a_{n}$ be a
sequence of real numbers such that $a_{n}\rightarrow 0$ and
$na_{n}\rightarrow \infty $. Assume that $\Vert D_{j,n}\Vert
_{\infty } = o(\sqrt{na_{n}})$ and that for all $\delta >0$, there
exists $ \sigma^2 \geq 0$ such that
\begin{equation}
\limsup_{n\rightarrow \infty }a_{n}\mathrm{log}\,\mathbf{P}\left( \Big|
\frac{1}{n}\sum_{j=1}^{k_{n}}\mathbf{E}(D_{j,n}^{2}|\mathcal{F}_{(j-1),n})- \sigma^2
\Big|\geq \delta \right) =-\infty \,.  \label{L2}
\end{equation}
Then, for the given sequence $a_n$ the partial sums processes
$Z_{n}(.)$ satisfy (\ref{mdpkey})
with the good rate function $I_\sigma(\cdot)$ defined in (\ref{rate}).
\end{lemma}

To be able to obtain the moderate deviation principle by approximation with
martingales we state next a simple approximation
lemma from Dembo and Zeitouni (1998, Theorem 4.2.13. p 130),
called exponentially
equivalence lemma.

\begin{lemma}\label{approlm}
Let $\xi _{n}(.):=\{\xi_n(t) \,, t \in [0,1] \}$ and
$\zeta _{n}(.):=\{\zeta_n(t) \, , t \in [0,1] \}$ be two processes
in $D([0,1])$. Assume that for any $\delta > 0$,
\begin{equation*}
\limsup_{n\rightarrow \infty }a_{n}\mathrm{log}\,{\mathbf P}(\sqrt{a_n}\sup_{t \in
[0,1]} | \xi _{n}(t)-\zeta _{n}(t)| \geq \delta )=-\infty
\end{equation*}
Then, if the sequence of processes $\xi _{n}(.)$ satisfies (\ref{mdpkey}) then so does
the sequence of processes $\zeta_{n}(.)$.
\end{lemma}

In dealing with dependent random variables, to brake
the dependence, a standard procedure is to divide first the
variables in blocks. This technique introduces a new parameter,
and so, in order to use a blocking procedure followed by a martingale
approximation, we have to establish a more specific exponentially
equivalent approximation, as stated in the following lemma:

\begin{lemma}
\label{simple} For any positive integer $m$, let $k_{n,m}$ be an
increasing sequence of integers going to infinity. Let
$\{d_{j,n}^{(m)}\}_{1\leq j\leq k_{n,m}}$ be a sequence of
triangular array of martingale differences adapted to a filtration
${\mathcal{F}^{(m)}_{j,n}}$. Define the normalized partial sums
process
$Z_{n}^{(m)}(t)=n^{-1/2}\sum_{i=1}^{[k_{n,m}t]}d_{i,n}^{(m)}$. Let
$a_{n}$ be a sequence of positive numbers such that
$a_{n}\rightarrow 0$ and $na_{n}\rightarrow \infty $. Assume that
for all $m\geq 1$
\begin{equation}
\sup_{1\leq j\leq k_{n,m}}\Vert d_{j,n}^{(m)}\Vert _{\infty }=o(\sqrt{na_{n}}
)\text{ as }n\rightarrow \infty  \label{Bound2}
\end{equation}
and that for all $\delta >0$, there exists $\sigma^2 \geq 0$ such
that
\begin{equation}
\lim_{m\rightarrow \infty }\lim \sup_{n\rightarrow \infty }\ a_{n}
\mathrm{log}\,\mathbf{P}\left( \Big|\frac{1}{n}\sum_{j=1}^{k_{n,m}}
\mathbf{E}((d_{j,n}^{(m)})^{2}|\mathcal{F}^{(m)}_{(j-1),n})-\sigma^2\Big|\geq
\delta \right) =-\infty \,.  \label{c11}
\end{equation}
Let $\{\zeta_{n}(t) \, , t \in [0,1]\}$ be a
sequence of $D[0,1]$--valued random variables such that for all
$\delta > 0$,
\begin{equation}
\lim_{m\rightarrow \infty }\lim \sup_{n\rightarrow \infty
}a_{n}\mathrm{log}\,{\mathbf P}(\sqrt{a_n}\sup_{t \in [0,1]}\vert \zeta _{n}(t) -
Z_{n}^{(m)}(t)\vert \geq \delta )=-\infty \label{close}
\end{equation}
Then, the processes $\zeta_{n}(.)$ satisfy (\ref{mdpkey})
with the good rate function $I_\sigma(\cdot)$ defined in (\ref{rate}).
\end{lemma}
{\bf Proof}. Define the functions
\begin{align*}
A_{1}(\delta ,n,m)& =a_{n}\mathrm{log}\,\mathbf{P} (\sup_{t \in [0,1]}\vert
\zeta _{n}(t)-Z_{n}^{(m)}(t)\vert \geq \delta )\;;\; \\
A_{2}(\delta ,n,m)& =a_{n}\mathrm{log}\, \mathbf{P}\left( \sup_{t \in [0,1] } \Big|\frac{1}{n}
\sum_{j=1}^{[k_{n,m}t]}\mathbf{E}((d_{j,n}^{(m)})^{2}
|\mathcal{F}^{(m)}_{(j-1),n})-t\sigma^2 \Big|\geq \delta \right) \\
A_{3}(\delta ,n,m)& =\mathrm{log}\,(\sup_{1\leq j\leq k_{n,m}}\Vert d_{j,n}^{(m)}
\Vert_{\infty })-\mathrm{log}\,(\sqrt{a_{n}n})\;.
\end{align*}
Observe that the functions $A_i, i=1,2,3$ satisfy the conditions
of Lemma \ref{convergence} from Appendix and so, we can find a sequence $m_{n}\rightarrow
\infty $ such that the martingale difference sequence
$(d_{j,n}^{(m_n)})$ satisfies the conditions of Lemma
\ref{triangmart}. We then derive that the sequence of processes
$Z_{n}^{(m_{n})}(.)$ satisfies (\ref{mdpkey}) and, by applying
Lemma \ref{approlm}, so does the sequence $\zeta _{n}(.)$ . \quad
$\diamond$
\subsection{Proof of Theorem \ref{mainresult}}
$\quad \;$
Let $m$ be an integer and $k=k_{n,m}=[n/m]$ (where, as before, $[x]$ denotes
the integer part of $x$).

The initial step of the proof is to divide the
variables in blocks of size $m$ and to make the sums in each block
\begin{equation*}
X_{i,m}=\sum_{j=(i-1)m+1}^{im}X_{j}\,,\text{ }i\geq 1.
\end{equation*}

Then we construct the martingales,
\begin{equation*}
M_{k}^{(m)}=\sum_{i=1}^{[n/m]}(X_{i,m}-\mathbf{E}(X_{i,m}|
\mathcal{F}_{(i-1)m}):=\sum_{i=1}^{[n/m]}D_{i,m}\,
\end{equation*}
and we define the process $\{M_{k}^{(m)}(t):\,t\in \lbrack 0,1]\}$ by
\begin{equation*}
M_{k}^{(m)}(t):=M^{(m)}_{[kt]}\,.
\end{equation*}

Now, we shall use Lemma \ref{simple} applied with $d_{j,n}^{(m)}=D_{j,m},$ and verify
the conditions (\ref{c11}) and (\ref{close}).

We start by proving (\ref{c11}). Notice first that
$\{D_{i,m}\}_{i\geq 1}$ is a rowwise stationary sequence of bounded
martingale differences. We have to verify
\begin{equation}
\lim_{m\rightarrow \infty }\lim \sup_{n\rightarrow \infty }a_{n}
{\rm log}\,\mathbf{P}\left( \Big |\frac{1}{n}\sum_{j=1}^{[n/m]}\mathbf{E}(D_{j,m}^{2}|
\mathcal{F}_{(j-1)m})- \sigma ^{2} \Big |\geq \delta \right)
=-\infty \,. \label{3}
\end{equation}
Notice that
\begin{equation*}
\mathbf{E}(D_{j,m}^{2}|\mathcal{F}_{(j-1)m})=\mathbf{E}(X_{j,m}^{2}|
\mathcal{F}_{(j-1)m})-(\mathbf{E}(X_{j,m}|\mathcal{F}_{(j-1)m)}))^{2}
\end{equation*}
and that, by stationarity
\begin{equation*}
\frac{1}{n}\Big \| \sum_{j=1}^{[n/m]}(\mathbf{E}(X_{j,m}|
\mathcal{F}_{(j-1)m}))^{2}\Big \| _{\infty }\leq \frac{\Vert \mathbf{E}(S_{m}|
\mathcal{F}_{0})\Vert _{\infty }^{2}}{m}\, .
\end{equation*}
Also
\begin{equation*}
 \Big \| \frac{1}{n}\sum_{j=1}^{[n/m]}\mathbf{E}(X_{j,m}^{2}|
\mathcal{F}_{(j-1)m})-\sigma ^{2}\Big \|_{\infty } \leq
\|m^{-1}\mathbf{E }(S_{m}^{2}|\mathcal{F}_{0})-\sigma ^{2}\|_{\infty
} +(1-km/n)\sigma ^{2} .
\end{equation*}
 Consequently
\begin{equation*}
\limsup_{n \rightarrow \infty} \Big \|\frac{1}{n}\sum_{j=1}^{[n/m]}
(\mathbf{E} (D_{j,m}^{2}|\mathcal{F}_{(j-1)m})-\sigma ^{2}\Big
\|_{\infty }\leq  \frac{\Vert \mathbf{E}(S_{m}|
\mathcal{F}_{0})\Vert _{\infty }^{2}}{m} + \|m^{-1}\mathbf{E
}(S_{m}^{2}|\mathcal{F}_{0})-\sigma ^{2}\|_{\infty }
\end{equation*}
which is smaller than $\delta/2$ provided $m$ is large enough,  by
the first part of Lemma \ref{techlemma} from Appendix and condition
(\ref{s2inf}). This proves (\ref{3}).

It remains to prove (\ref{close}), that means in our notation that
for any $\delta >0$
\begin{equation}
\lim_{m\rightarrow \infty }\lim \sup_{n\rightarrow \infty }a_{n}
{\rm log}\,\mathbf{P}\left( \sqrt{\frac{a_{n}}{n}}\sup_{t\in \lbrack
0,1]}|S_{[nt]}-M_{k}^{(m)}(t)|\geq \delta \right) =-\infty \,.
\label{1}
\end{equation}
Notice first that
\begin{align*}
\sup_{t\in \lbrack 0,1]}|S_{[nt]}-M_{k}^{(m)}(t)|& \leq \sup_{t\in
\lbrack 0,1]}|\sum_{i=[k \, t]m+1}^{[nt]}X_{i}|+\sup_{t\in \lbrack
0,1]}|\sum_{i=1}^{[k \, t]}\mathbf{E}(X_{i,m}|\mathcal{F}_{(i-1)m})| \\
& \leq o(\sqrt{na_{n}})+\max_{1\leq j\leq
[n/m]}|\sum_{i=1}^{j}\mathbf{E} (X_{i,m}|\mathcal{F}_{(i-1)m})|.
\end{align*}
Then, by using Lemma \ref{puw} we derive that
\begin{align*}
& a_{n}\mathrm{log}\,\mathbf{P}\left(
\sqrt{\frac{a_{n}}{n}}\max_{1\leq j\leq
[n/m]}|\sum_{i=1}^{j}\mathbf{E}(X_{i,m}|\mathcal{F}_{(i-1)m})|\geq
\delta \right) \\
& \quad \leq a_{n}\mathrm{log}(4\sqrt{e})-
\frac{\delta ^{2}m}{2(\Vert \mathbf{E}(S_{m}|\mathcal{F}_{0})\Vert _{\infty
}+80\sum_{j=1}^{\infty }j^{-3/2}\Vert \mathbf{E}(S_{jm}|\mathcal{F}
_{0})\Vert _{\infty })^{2}}\,.
\end{align*}
which is convergent to $-\infty $ when $n\rightarrow \infty $ followed by
$m\rightarrow \infty $, by Lemma \ref{techlemma}.

\subsection{Proof of Corollary  \ref{cor1} and Remark 5}
$\quad \;$
Notice that obviously, by triangle inequality
and changing the order of summation, (\ref{bis}) implies
(\ref{mw}). So, in order to establish
both Corollary  \ref{cor1} and Remark 5, we just have to
show that condition (\ref{mw}) together with (\ref{mix}) imply condition (\ref{s2inf}).
This will be achieved by using the following two lemmas.

First let us introduce some notations.
Let $S_{a,b} =S_{b}-S_{a}$ and set
\begin{eqnarray*}
\widetilde \Delta _{r, \infty} = \sum_{j=r}^{\infty}2^{-j/2}\Vert
\mathbf{E}(S_{2^{j}}|\mathcal{F} _{0})\Vert _{\infty }\;,\;
\Delta_{\infty} =\Vert E(X_1^2|\mathcal{F}_0)\Vert_{\infty}^{1/2} +
\sum_{j=0}^{\infty }2^{-j/2}\Vert \mathbf{E}
(S_{2^{j}}|\mathcal{F}_{0})\Vert _{\infty }\;.
\end{eqnarray*}
By Peligrad and Utev (2005), $\widetilde \Delta _{0, \infty}<\infty $ is
equivalent to (\ref{mw}).
\begin{lemma}\label{PUvar}
\label{infbound} Assume that $X_0$ is ${\mathcal F}_0$--measurable
and that $\Vert \mathbf E(X_{1}^2|\mathcal{F}_{0})\Vert _{\infty
}<\infty $. Let $n$, $r$ be integers such that $2^{r-1}<n\leq
2^{r}$. Then
\begin{equation*}
\|\mathbf{E}(S_{n}^{2}|\mathcal{F}_{0})\|_{\infty }\leq n\Big(\Vert
E(X_{1}^2|\mathcal{F}_{0})\Vert _{\infty }^{1/2}+ \tfrac{1}{2}
\sum_{j=0}^{r-1}2^{-j/2}\Vert \mathbf{E}(S_{2^{j}}|\mathcal{F}
_{0})\Vert _{\infty }\Big)^{2} \leq n\Delta_{\infty} ^{2} \, .
\end{equation*}
Moreover, under (\ref{mw}),
\begin{equation*}
\Vert n^{-1}\mathbf{E}(S_{n}^{2}|\mathcal{I})-\eta \Vert _{\infty
}\rightarrow 0\quad \mbox{as}\quad n\rightarrow \infty \, ,
\end{equation*}
where ${\cal I}$ is the $\sigma$-field of all $T$-invariant sets and
\begin{equation*}
\eta =\mathbf{E}(X_{1}^{2}|\mathcal{I})+\sum_{j=0}^{\infty }2^{-j}\mathbf{E}
(S_{2^{j}}(S_{2^{j+1}}-S_{2^{j}})|\mathcal{I}) \, .
\end{equation*}
In particular, if ${\mathbf E}(X_i X_j |{\mathcal
F}_{-\infty})={\mathbf E}(X_i X_j)$, for any $i, j$ in ${\mathbf
Z}$, then
\begin{equation*}
\eta=\sigma ^{2}=\mathbf{E}(X_{1}^{2})+\sum_{j=0}^{\infty }2^{-j}\mathbf{E}
(S_{2^{j}}(S_{2^{j+1}}-S_{2^{j}}))\, .
\end{equation*}
\end{lemma}
{\bf Proof.} The proofs of the first three statements are almost
identical to the proof of the corresponding facts in Proposition
2.1 of Peligrad and Utev (2005). The only changes are to replace
everywhere the $L_{2}$--norm $\Vert x\Vert $ by the $L_{\infty
}$--norm $\Vert x\Vert _{\infty }$ and the usual expectation
$\mathbf{E(}X)$ by the conditional expectation $\mathbf{E(}X) =
\mathbf{E}(X|\mathcal{F}_{0})$. The last statement follows from
Proposition 2.12 in Bradley (2002), since for all $i,j$,
\begin{equation*}
\mathbf{E}(X_{i}X_{j}|\mathcal{I})=\mathbf{E}(\mathbf{E}(X_{i}X_{j}|
\mathcal{F}_{-\infty })| \mathcal{I})=\mathbf{E}(X_{i}X_{j})\, .
\quad \diamond
\end{equation*}
\begin{lemma}
\label{inferg2} Assume that $X_0$ is ${\mathcal F}_0$--measurable
and that $\Vert \mathbf{E}(X_{1}^2|\mathcal{F}_{0})\Vert _{\infty
}<\infty$. Suppose that the conditions (\ref{mw}) and (\ref{mix})
are satisfied. Then,
\begin{equation*} \Vert
n^{-1}\mathbf{E}(S_{n}^{2}|\mathcal{F}_{0})-\sigma^2 \Vert _{\infty
}\rightarrow 0\quad \mbox{as}\quad n\rightarrow \infty
\end{equation*}
\end{lemma}
{\bf Proof}. By Lemma \ref{infbound},
it is enough to show that
\begin{eqnarray*}
\frac{1}{n}\Vert\mathbf{E}(S_{n}^{2}|\mathcal{F}_{0})-\mathbf{E}(S_{n}
^{2})\Vert_{\infty}\rightarrow0\text{ as }n\rightarrow\infty
\end{eqnarray*}
We prove this lemma by diadic recurrence. For $t$ integer, denote
\begin{eqnarray*}
A_{t,k}=\Vert\mathbf{E}(S_{t}^{2}|\mathcal{F}_{-k})-\mathbf{E}(S_{t}^{2}
)\Vert_{\infty} \, .
\end{eqnarray*}
Then, by the properties of conditional expectation and stationarity, for all
$t\geq1$
\begin{eqnarray*}
\ A_{2t,k} &=&\Vert\mathbf{E}(S_{2t}^{2}|\mathcal{F}_{-k})-\mathbf{E}
(S_{2t}^{2})\Vert_{\infty}\leq
\Vert\mathbf{E}(S_{t}^{2}|\mathcal{F}_{-k})-\mathbf{E}(S_{t}^{2})
\Vert_{\infty}\\
&&\qquad + \Vert\mathbf{E}(S_{t,2t}^{2}|\mathcal{F}_{-k})-\mathbf{E}
(S_{t}^{2})\Vert_{\infty} +
2\Vert\mathbf{E}(S_{t}S_{t,2t}|\mathcal{F}_{-k})-\mathbf{E}(S_{t}
S_{t,2t})\Vert_{\infty}\\
&\leq&
2\Vert\mathbf{E}(S_{t}^{2}|\mathcal{F}_{-k})-\mathbf{E}(S_{t}^{2}
)\Vert_{\infty}+2\Vert\mathbf{E}(S_{t}S_{t,2t}|\mathcal{F}_{-k})\Vert_{\infty
}+2|\mathbf{E}(S_{t}S_{t,2t})| \, .
\end{eqnarray*}
Using for the last two terms the bound from Lemma \ref{infbound}, the
Cauchy-Schwartz inequality and stationarity, we have
\begin{eqnarray*}
A_{2t,k}\leq2A_{t,k}+4t^{1/2}\Delta_{\infty}\Vert\mathbf{E}(S_{t}
|\mathcal{F}_{0})\Vert_{\infty}
\end{eqnarray*}
Whence, with the notation
\begin{eqnarray*}
B_{r,k}=2^{-r}\Vert\mathbf{E}(S_{2^{r}}^{2}|\mathcal{F}_{-k})-\mathbf{E}
(S_{2^{r}}^{2})\Vert_{\infty}=2^{-r}A_{2^{r},k}
\end{eqnarray*}
by recurrence, for all $r\geq m$ and all $k>0$, we derive
\begin{eqnarray*}
B_{r,k}  \leq B_{r-1,k}+2^{\frac{-r+3}{2}}\Delta_{\infty}\Vert
\mathbf{E}(S_{2^{r-1}}|\mathcal{F}_{0})\Vert_{\infty}
\leq B_{m,k}+2\Delta_{\infty}
\sum_{j=m}^{r}2^{-j/2}\Vert\mathbf{E}(S_{2^{j}
}|\mathcal{F}_{0})\Vert_{\infty}\,.
\end{eqnarray*}
Therefore
\begin{equation}
2^{-r}\Vert\mathbf{E}(S_{2^{r}}^{2}|\mathcal{F}_{-k})-\mathbf{E}(S_{2^{r}}
^{2})\Vert_{\infty}\leq B_{m,k}+2\Delta_{\infty}\widetilde{\Delta}_{m,\infty
}\,.\label{bounds2r}
\end{equation}
Now notice that, by stationarity and triangle inequality
\begin{eqnarray}
\Vert\mathbf{E}(S_{2^{r}}^{2}|\mathcal{F}_{0})-\mathbf{E}(S_{2^{r}}^{2}
)\Vert_{\infty}  \leq \Vert\mathbf{E}(S_{2^{r}}^{2}|\mathcal{F}
_{-k})-\mathbf{E}(S_{2^{r}}^{2})\Vert_{\infty}
+\Vert\mathbf{E}(S_{2^{r}}^{2}-S_{k,k+2^{r}}^{2}|\mathcal{F}_{0}
)\Vert_{\infty}\,,\label{bounds3r*}
\end{eqnarray}
and that by Lemma \ref{infbound}
\begin{eqnarray} \label{bounds2r**}
\Vert\mathbf{E}(S_{2^{r}}^{2}-S_{k,k+2^{r}}^{2}|\mathcal{F}_{0})\Vert
_{\infty}
&\leq&\Vert\mathbf{E}((S_{2^{r}}-S_{k,k+2^{r}})^{2}|\mathcal{F}
_{0})\Vert_{\infty}^{1/2}
\Vert\mathbf{E}((S_{2^{r}}+S_{k,k+2^{r}})^{2}|\mathcal{F}_{0})\Vert_{\infty
}^{1/2}\nonumber \\
&\leq&
4k\Vert\mathbf{E}(X_{1}^{2}|\mathcal{F}_{0})\Vert_{\infty}^{1/2}
\Vert\mathbf{E}(S_{2^{r}}^{2}|\mathcal{F}_{0})\Vert_{\infty}^{1/2}
\nonumber \\ & \leq &
2^{2+r/2}k\Vert\mathbf{E}(X_{1}^{2}|\mathcal{F}_{0})\Vert_{\infty}^{1/2}
\Delta_{\infty}\,.
\end{eqnarray}
Then, starting from (\ref{bounds3r*}) and using (\ref{bounds2r}) and
(\ref{bounds2r**}), we derive that for $r\geq m+1$,
\begin{eqnarray*}
2^{-r}\Vert\mathbf{E}(S_{2^{r}}^{2}|\mathcal{F}_{0})-\mathbf{E}(S_{2^{r}}
^{2})\Vert_{\infty}\leq B_{m,
k}+2\Delta_{\infty}\widetilde{\Delta}_{m,\infty }+2^{-r/2+2}k\Vert
\mathbf{E}(X_{1}^2|\mathcal{F}_{0})\Vert _{\infty
}^{1/2}\Delta_{\infty}\,.
\end{eqnarray*}
As a consequence
\begin{eqnarray*}
\limsup_{r\rightarrow\infty}2^{-r}\Vert\mathbf{E}(S_{2^{r}}^{2}|\mathcal{F}
_{0})-\mathbf{E}(S_{2^{r}}^{2})\Vert_{\infty}\leq B_{m,
k}+2\Delta_{\infty }\widetilde{\Delta}_{m,\infty}\,.
\end{eqnarray*}
Then, we first let $k\rightarrow\infty$ and by Condition \ref{mix} it follows
that $\lim_{k\rightarrow\infty}B_{m,k}=0$. Then, we let $m$ tend to infinity
and by Condition (\ref{mw}), we derive
\begin{eqnarray*}
\lim_{r\rightarrow\infty}2^{-r}\Vert\mathbf{E}(S_{2^{r}}^{2}|\mathcal{F}
_{0})-\mathbf{E}(S_{2^{r}}^{2})\Vert_{\infty}=0\,.
\end{eqnarray*}
To complete the proof of the lemma we use the
diadic expansion $n=\Sigma_{k=0}^{r-1}2^{k}a_{k}$ where $a_{r-1}=1$
and $a_{k}\in\{0,1\}$ and continue the proof as in Proposition
2.1 in Peligrad and Utev (2005).\quad
$\diamond$

\subsection{Proof of Theorem \ref{projresult}}
$\quad \;$
Fix a positive integer $m$ and define the stationary sequence
\begin{equation*}
\xi
_{j,m}:=\mathbf{E}(X_{j}|\mathcal{F}_{j+m-1})-\mathbf{E}(X_{j}|\mathcal{F}_{j-m})
\end{equation*}
Using a standard martingale decomposition (see also Hall and Heyde, 1980), we define
\begin{equation*}
\theta _{j,m}=\sum_{t=0}^{\infty }\mathbf{E}(\xi _{j+t,m}|\mathcal{F}_{j+m-1})=
\sum_{k=0}^{2m-2}\mathbf{E}(\xi _{j+k,m}|\mathcal{F}_{j+m-1}) \, .
\end{equation*}
and observe that
\begin{eqnarray}\label{newadd1}
\Vert \theta _{0,m}\Vert
_{\infty }= \Vert \sum_{k=0}^{2m-2}\sum_{i=k-m +
1}^{m-1}P_i(X_k)\Vert_{\infty}\leq 2m\sum_{i \in {\mathbf Z}} \Vert
P_0(X_i) \Vert_{\infty}<\infty\,.
\end{eqnarray}
Then, $\mathbf{E}(\theta _{j+1,m}|\mathcal{F}_{j+m-1}) =\theta _{j,m}-\xi _{j,m} $
and thus,
\begin{equation} \label{deccobord}
\sum_{j=1}^{k}\xi _{j,m}=\theta _{1,m}-\theta
_{k+1,m}+\sum_{j=1}^{k}d_{j,m} \, .
\end{equation}
where $d_{j,m}:=\theta _{j+1,m}-\mathbf{E}(\theta
_{j+1,m}|\mathcal{F}_{j+m-1})$ is a stationary bounded
martingale difference.

Moreover,
\begin{equation} \label{decSn}
\sum_{j=1}^{k}X_j = \sum_{j=1}^{k}d_{j,m} + R_{k,m}\, ,
\end{equation}
where $$ R_{k,m}:= \theta _{1,m}-\theta _{k+1,m} +
\sum_{j=1}^{k}[X_j - E(X_{j}|\mathcal{F}_{j+m-1})+
\mathbf{E}(X_{j}|\mathcal{F}_{j-m})] \, .$$ First, we show that
$R_{k, m}$ is negligible for the moderate deviation. We notice that by (\ref{newadd1}) it is enough to establish that $$ R'_{k,m}:=
\sum_{j=1}^{k}[X_j - E(X_{j}|\mathcal{F}_{j+m-1})+
\mathbf{E}(X_{j}|\mathcal{F}_{j-m})] \, $$ is negligible.
Observe that
\begin{eqnarray}
&&X_j - E(X_{j}|\mathcal{F}_{j+m-1})+
\mathbf{E}(X_{j}|\mathcal{F}_{j-m}) = \sum_{|t|\geq m} P_{j-t}(X_j)
\quad\mbox{and}\nonumber\\
&&\sum_{j\in {\mathbf {Z}}}\Vert P_0 \big ( X_j -
E(X_{j}|\mathcal{F}_{j+m-1})+ \mathbf{E}(X_{j}|\mathcal{F}_{j-m})
\big ) \Vert_{\infty} \leq \sum_{|k| \geq m}^{\infty }\Vert
P_{0}(X_{k})\Vert _{\infty } =:D_m \label{newadd2}
\end{eqnarray}
Now, the exponential inequality given in
Lemma \ref{exp2} entails that
\begin{equation*}
{\mathbf P}\Big(\max_{1\leq k\leq n}\Big|
\sum_{j=1}^{k}\mathbf{E}
[X_j - E(X_{j}|\mathcal{F}_{j+m-1})+
\mathbf{E}(X_{j}|\mathcal{F}_{j-m})]
\Big|\geq \delta \sqrt{n/a_{n}}\Big)\leq 8\exp \Big(-\frac{\delta
^{2}n}{a_{n}2nD_{m}^{2}} \Big)
\end{equation*}
The last inequality together with (\ref{projcond}) and Lemma \ref{approlm}
reduces the theorem to the MDP principle for
bounded stationary martingale difference
$\{d_{j,m}\;;\;j\in \Z\}$.

Then, by Lemma \ref{simple}, it remains to  verify that
\begin{equation*}
\lim_{m\rightarrow \infty }\limsup_{n\rightarrow \infty }a_{n}\ln {\mathbf P}
\Big(\Big|\frac 1n \sum_{j=1}^{n}(\mathbf
E(d_{j,m}^{2}|\mathcal{F}_{j+m-1})-\sigma^{2})\Big|\geq \delta
\Big)=-\infty\,.
\end{equation*}
In order to prove this convergence, by Lemma \ref{techlemma1}, applied with
$\ B = 2 \big ( \sum_{\ell \in {\mathbf Z}} \Vert P_{0}(X_{\ell})
\Vert_{\infty} \big) ^2$, it is enough to establish that
\begin{equation*}
\lim_{m\rightarrow \infty }\limsup_{n\rightarrow \infty }
\Big\|\frac{1}{n}\sum_{j=1}^{n}
(\mathbf{E}(d_{j,m}^{2}|\mathcal{F}_{m-1})- \sigma^2 )
\Big\|_{\infty }= 0 \, .
\end{equation*}
Since $\{d_{j,m}\}$ is a martingale difference, it follows from
the decomposition (\ref{deccobord}) and (\ref{newadd1}), that
it remains to prove that
\begin{equation} \label{doublelim}
\lim_{m\rightarrow \infty }\limsup_{n\rightarrow \infty }
\Big\|\frac{1}{n} {\mathbf E} \Big ( \Big (
\sum_{j=2m-1}^{n}\xi_{j,m} \Big )^{2}|\mathcal{F}_{m-1} \Big )-
\sigma^2 \Big\|_{\infty }= 0 \, .
\end{equation}
Write
$$ \Big ( \sum_{j=2m-1}^{n}\xi_{j,m} \Big )^{2} =
\sum_{i=2m-1}^{n}\xi_{i,m}^2 + 2 \sum_{i=2m-1}^n \sum_{j =
i+1}^{(N+i)\wedge n} \xi_{i,m}\xi_{j,m} + 2 \sum_{i=2m-1}^n \sum_{j=
N+i+1}^n \xi_{i,m}\xi_{j,m} \, .$$
Notice that, since $\xi_{j,m} = \Sigma_{k=
j-m+1}^{j+m-1}P_k(X_j)$, we get
\begin{eqnarray*}
& &\frac {1}{n}
\Vert \sum_{i=2m-1}^n \sum_{j= N+i+1}^n  {\mathbf E} \big (
\xi_{i,m}\xi_{j,m}
 |\mathcal{F}_{m-1} \big )\Vert_{\infty} \leq
 \frac{1}{n}
\sum_{i=2m-1}^n \sum_{j= N+i+1}^{\infty}
\Vert {\mathbf E} \big (
\xi_{i,m}\xi_{j,m} |\mathcal{F}_{m-1} \big )\Vert_{\infty}
\\
& &
\leq  \frac{1}{n}
\sum_{i=2m-1}^n \sum_{k= i-m+1}^{i+m-1} \Vert P_k(X_i)
\Vert_{\infty}\sum_{\ell \geq N } \Vert P_k (X_{i+ \ell})
\Vert_{\infty} \;\leq\; \sum_{i \in {\mathbf Z}} \Vert P_0(X_i)
\Vert_{\infty} \sum_{|\ell |\geq N/2 } \Vert P_0 (X_{ \ell})
\Vert_{\infty} \, \to 0
\end{eqnarray*}
as $N\to \infty$, uniformly in $n$, and so, (\ref{doublelim}) is
implied by
\begin{equation} \label{doublelimbis}
\lim_{m\rightarrow \infty }\limsup_{n\rightarrow \infty }
\Big\|\frac{1}{n} {\mathbf E} \Big ( \sum_{i=2m-1}^{n}\xi_{i,m}^2 +
2 \sum_{i=2m-1}^n \sum_{j = i+1}^{(N+i)\wedge n} \xi_{i,m}\xi_{j,m}
|\mathcal{F}_{m-1} \Big )- \sigma_N^2 \Big\|_{\infty }= 0 \, ,
\end{equation}
where $\sigma^2_N = {\mathbf E} (X_0^2) +2{\mathbf E} (X_0X_1) +
\dots +2{\mathbf E} (X_0X_{N-1})$. Write $\xi_{i,m} = X_i +
(\xi_{i,m} - X_i)$. By condition (\ref{mix}), we easily get that
\begin{equation*} \label{doublelimter}
\lim_{m\rightarrow \infty }\limsup_{n\rightarrow \infty }
\Big\|\frac{1}{n} {\mathbf E} \Big ( \sum_{i=2m-1}^{n}X_i^2 + 2
\sum_{i=2m-1}^n \sum_{j = i+1}^{(N+i)\wedge n}
X_iX_j|\mathcal{F}_{m-1} \Big )- \sigma_N^2 \Big\|_{\infty }= 0 \, ,
\end{equation*}
hence (\ref{doublelimbis}) holds since
\begin{eqnarray*}
\Vert X_i -\xi_{i,m} \Vert_{\infty} \leq
\sum_{|k| \geq m}\Vert P_0 ( X_k) \Vert_{\infty}\to 0\quad
\mbox{as}\quad m\to\infty
\;.
\quad \diamond
\end{eqnarray*}

\subsection{Proof of Proposition \ref{proplip1}}
 Let ${\mathcal F}_k= \sigma(\varepsilon_i, i \leq k)$. From Theorem 4.4.7
in Berbee (1979), there exists $(\varepsilon'_i)_{i >0}$ distributed
as $(\varepsilon_i)_{i > 0}$ and independent of ${\mathcal F}_0$
such that
$$
 \|\E({\bf 1}_{\{\varepsilon_k\neq \varepsilon'_k, \,  \text{for some $k\geq
 n$}\}}|{\mathcal F}_0)\|_\infty=\phi_{\varepsilon}(n) \, .
$$
Let $(\varepsilon_i^{(0)})_{i \in \Z}$ be the sequence defined by
$\varepsilon_i^{(0)}=\varepsilon_i$ if $i \leq 0$ and
$\varepsilon_i^{(0)}=\varepsilon'_i$ if $i>0$. Let
$(\varepsilon_i^{(-1)})_{i \in \Z}$ be the sequence defined by
$\varepsilon_i^{(-1)}=\varepsilon_i$ if $i < 0$,
$\varepsilon_i^{(-1)}=\varepsilon'_i$ if $i>0$ and
$\varepsilon_0^{(-1)}=x$ where $x \in A$. Define now
$Z_k=H((\varepsilon_{k-i})_{i \in \Z})$,
$Z_k^{(0)}=H((\varepsilon^{(0)}_{k-i})_{i \in \Z})$ and
$Z_k^{(-1)}=H((\varepsilon^{(-1)}_{k-i})_{i \in \Z})$. We shall
apply Theorem \ref{projresult}. Note first that (\ref{projcond}) is
equivalent to $\sum_{i \in \Z} \|P_0(Z_i)\|_\infty < \infty$. Now
$$
P_0(Z_i)=\E(Z_i^{(0)}|{\cal F}_0)- \E(Z_i^{(-1)}|{\cal F}_{-1})
+\E(Z_i-Z_i^{(0)}|{\cal F}_0)- \E(Z_i-Z_i^{(-1)}|{\cal F}_{-1})\, .
$$
Denoting by ${\mathbf{E}}_{\varepsilon}(\cdot)$ the conditional
expectation with respect to $\varepsilon$, we infer from $C(A)$ that
$$
|\E(Z_i^{(0)}|{\cal F}_0)- \E(Z_i^{(-1)}|{\cal F}_{-1})|=
|{\mathbf{E}}_{\varepsilon}(H((\varepsilon^{(0)}_{i-j})_{j \in
\Z})-H((\varepsilon^{(-1)}_{i-j})_{j \in \Z}))|\leq \Delta_i \, .
$$
Now, from $C(A)$ again,
$$
|\E(Z_i-Z_i^{(0)}|{\cal F}_0)| \leq \sum_{k=1}^\infty \Delta_{i-k}
\E({\bf 1}_{\varepsilon_k \neq \varepsilon'_k}|{\cal F}_0) \
\text{and} \ \E(Z_i-Z_i^{(-1)}|{\cal F}_{-1})\leq \Delta_i +
\sum_{k=1}^\infty \Delta_{i-k} \E({\bf 1}_{\varepsilon_k \neq
\varepsilon'_k}|{\cal F}_{-1}).
$$
Consequently, by the $\phi$-mixing property, we obtain the upper
bound
$$
\sum_{i \in \Z} \|P_0(Z_i)\|_\infty \leq 2 \sum_{i \in \Z} \Delta_i
+2 \sum_{i \in \Z} \sum_{k=1}^\infty \Delta_{i-k} \phi_\varepsilon
(k)\, ,
$$
which is finite provided that $\sum_{i \in \Z} \Delta_i < \infty$
and $\sum_{k >0} \phi_\varepsilon (k) < \infty$. It remains to prove
(\ref{mix}). Let $X_k^{(0)}=Z_k^{(0)}-\E(Z_k^{(0)})$. We have
\begin{eqnarray}\label{debildec}
\Vert \mathbf{E} (X_kX_l |\mathcal{F}_{0})-\mathbf{E} ( X_kX_l )
\Vert_\infty&\leq& \Vert \mathbf{E} (X_k^{(0)}X_l^{(0)}
|\mathcal{F}_{0})-\mathbf{E} ( X_k^{(0)}X_l^{(0)} ) \Vert_\infty\nonumber\\
&+&\Vert \mathbf{E} (X_k(X_l-X_l^{(0)}) |\mathcal{F}_{0})-\mathbf{E}
( X_k(X_l-X_l^{(0)}) ) \Vert_\infty\nonumber\\
&+&\Vert \mathbf{E} (X_l^{(0)}(X_k-X_k^{(0)})
|\mathcal{F}_{0})-\mathbf{E} ( X_l^{(0)}(X_k-X_k^{(0)}) )
\Vert_\infty\, .
\end{eqnarray}
Clearly, by $C(A)$ and the $\phi$-mixing property,
$$
\Vert \mathbf{E} (X_k(X_l-X_l^{(0)}) |\mathcal{F}_{0})-\mathbf{E} (
X_k(X_l-X_l^{(0)}) ) \Vert_\infty \leq 4\|X_k\|_\infty
\sum_{k=1}^\infty \Delta_{l-k} \phi_\varepsilon (k)\, ,
$$
which tends to zero as $l$ tends to infinity. In the same way
$$
 \lim_{k \rightarrow 0} \Vert \mathbf{E} (X_l^{(0)}(X_k-X_k^{(0)})
|\mathcal{F}_{0})-\mathbf{E} ( X_l^{(0)}(X_k-X_k^{(0)}) )
\Vert_\infty=0 \, .
$$
Let $H_k=H-\E(Z_k^{(0)})$. Let $(\eta_i)_{i \in \Z}$ be distributed
as $(\varepsilon_i)_{i \in \Z}$ and independent of
$((\varepsilon_i)_{i \in \Z}, (\varepsilon'_i)_{i>0})$, and let
$(\eta_i^{(0)})_{i \in \Z}$ be the sequence defined by
$\eta_i^{(0)}=\eta_i$ if $i \leq 0$ and
$\eta_i^{(0)}=\varepsilon'_i$ if $i>0$. With this notations, we have
\begin{eqnarray}\label{debildec2}
 \mathbf{E} (X_k^{(0)}X_l^{(0)} |\mathcal{F}_{0})-\mathbf{E} (
X_k^{(0)}X_l^{(0)} )
&=&{\mathbf{E}}_{\varepsilon}(H_k((\varepsilon^{(0)}_{k-i})_{i \in
Z})(H_l((\varepsilon^{(0)}_{l-i})_{i \in
Z})-H_l((\eta^{(0)}_{l-i})_{i \in Z})))\, \nonumber \\
&+& {\mathbf{E}}_{\varepsilon}(H_l((\eta^{(0)}_{l-i})_{i \in
Z})(H_k((\varepsilon^{(0)}_{k-i})_{i \in
Z})-H_k((\eta^{(0)}_{k-i})_{i \in Z})))
\end{eqnarray}
Consequently, applying $C(A)$ once more, we have that
$$
\|\mathbf{E} (X_k^{(0)}X_l^{(0)} |\mathcal{F}_{0})-\mathbf{E} (
X_k^{(0)}X_l^{(0)} )\|_\infty \leq \|X_k^{(0)}\|_\infty \sum_{i \geq
l} \Delta_i + \|X_l^{(0)}\|_\infty \sum_{i \geq k} \Delta_i \, ,
$$
which tends to zero as $k$ and $l$ tends to infinity. This completes
the proof. \quad $\diamond $

\subsection{Proof of Proposition \ref{proplip2}}
We shall apply Corollary \ref{cor1}.
We use the same notations as for the proof of  Proposition
\ref{proplip1}. With these notations, we have
$$
\E(X_k |{\cal F}_0)= \E(X_k^{(0)} |{\cal F}_0) + \E(X_k-X_k^{(0)}
|{\cal F}_0)\, .
$$
Now, applying $C'(A)$,
$$
|\E(X_k-X_k^{(0)} |{\cal F}_0)|\leq \sum_{i=1}^k R_{k-i}\E({\bf
1}_{\varepsilon_i\neq \varepsilon'_i}|{\mathcal F}_0) + \sum_{i=1}^k
R_{k-i}{\bf P}(\varepsilon_i\neq \varepsilon'_i)\, ,
$$
and by the $\phi$-mixing property,
\begin{equation}\label{ineq40}
\|\E(X_k-X_k^{(0)} |{\cal F}_0)\|_\infty \leq 2\sum_{i=1}^k
R_{k-i}\phi_\varepsilon(i) \, .
\end{equation}
Now, by $C'(A)$ again,
$$
\|\E(X_k^{(0)} |{\cal
F}_0)\|_\infty=\|\E_\varepsilon(H((\varepsilon^{(0)}_{k-i})_{i \in
\Z})-H((\eta^{(0)}_{k-i})_{i \in \Z}))\|_\infty \leq R_k \, .
$$
Consequently, since $\phi_\varepsilon (0)>0$, the condition
(\ref{bis}) is implied by (\ref{mixrphi}). It remains to prove
(\ref{mix}). We start from  the decomposition (\ref{debildec}). By  (\ref{ineq40}),
\begin{eqnarray*}
\Vert \mathbf{E} (X_k(X_l-X_l^{(0)}) |\mathcal{F}_{0})-\mathbf{E} (
X_k(X_l-X_l^{(0)}) ) \Vert_\infty &\leq & 4\|X_k\|_\infty
\sum_{i=1}^l
R_{l-i}\phi_\varepsilon(i)\,, \\
\text{and}\quad \Vert \mathbf{E} (X_l^{(0)}(X_k-X_k^{(0)})
|\mathcal{F}_{0})-\mathbf{E} ( X_l^{(0)}(X_k-X_k^{(0)}) )
\Vert_\infty &\leq & 4\|X_l^{(0)}\|_\infty
\sum_{i=1}^k
R_{k-i}\phi_\varepsilon(i)\, .
\end{eqnarray*}
Hence, in view of (\ref{mixrphi}), these two terms converges to zero
as $k$ and $l$ tend to infinity. From (\ref{debildec2}) and condition $C'(A)$, we have that
$$
\|\mathbf{E} (X_k^{(0)}X_l^{(0)} |\mathcal{F}_{0})-\mathbf{E} (
X_k^{(0)}X_l^{(0)} )\|_\infty \leq \|X_k^{(0)}\|_\infty R_l + \|X_l^{(0)}\|_\infty R_k\, ,
$$
which again converges to zero as $k$ and $l$ tend to infinity. This completes the proof. \quad $\diamond $

\subsection{Proof of Proposition \ref{contracting}}
It suffices to prove that for any $f$ in ${\cal L}$, the sequence
$X_i= f(Y_i)-\mu(f)$ satisfies the conditions (\ref{bis}) and
(\ref{mix}) of Corollary \ref{cor1}.

Note first that (\ref{mix}) holds because of (\ref{contrac2}) and
because any continuous function from $[0,1]$ to ${\mathbf R}$ can be
uniformly approximated by Lipschitz functions.

From  Lemma \ref{rhoarith}, we have that
$$
\|K^n(f)-\mu(f)\|_{\infty, \mu} \leq c(C \rho^n) \, ,
$$
for some concave non decreasing function $c$. Consequently
(\ref{bis}) holds as soon as $\sum_{k>0} k^{-1/2} c (C\rho^k)$ is
finite, which in turn is equivalent to (\ref{intc}).

\subsection{Proof of Lemma \ref{rhoarith}}
$\quad \;$ Let $(Y_i)_{i \geq 1}$ be the Markov chain with
transition Kernel $K$ and and invariant measure $\mu$.
 From Lemma 1 in Dedecker
and Merlev\`ede (2006), we know that there exists $Y_k^*$
distributed as $Y_k$ and independent of $Y_0$ such that
$$
      \sup_{g \in \Lambda_1} \| K^k(g)-\mu(g)  \|_{\infty, \mu} =
    \|{\mathbf E}(|Y_k-Y_k^*||Y_0)\|_\infty \, .
$$
For any $f$ such that $|f(x) - f(y) | \leq c ( |x-y|)$, we have
\begin{eqnarray*}
 \| K^k(f)-\mu(f)\|_{\infty, \mu} &= & \|{\mathbf E}(f(Y_k)|Y_0) - {\mathbf E}(f(Y_k^*)|Y_0) \|_\infty \\
 &\leq &    \|{\mathbf E}(c(|Y_k-Y_k^*|)|Y_0)\|_\infty \, .
 \end{eqnarray*}
 Since $c$ is concave and non decreasing, we get that
 \begin{eqnarray*}
 \| K^k(f)-\mu(f)\|_{\infty, \mu}
 \leq    \|c ( {\mathbf E}(|Y_k-Y_k^*||Y_0) ) \|_\infty  \leq
 c ( \|{\mathbf E}(|Y_k-Y_k^*||Y_0)\|_\infty ) \,
 ,
 \end{eqnarray*}
and the proof is complete. \quad $\diamond $

\subsection{Proof of Corollary \ref{proexp}}
$\quad \;$ Let $(Y_i)_{i \geq 1}$ be the Markov chain with
transition Kernel $K$ and invariant measure $\mu$. Using the
equation (\ref{eq}) is easy to see that $(Y_0, \ldots, Y_n)$ is
distributed as $(T^{n+1}, \ldots, T)$. Consequently, for $f$ in
${\cal L}$, Corollary \ref{proexp} follows from Proposition
\ref{contracting} and Condition (\ref{expan}).

Assume now that $f$ is $BV$. We shall prove that the sequence $X_i=
f(Y_i)-\mu(f)$ satisfies the conditions (\ref{bis}) and (\ref{mix})
of Corollary \ref{cor1}. Since $K$ is $BV$-contracting we have that
$$
  \| {\mathbf E}(X_k|Y_0) \|_{\infty} = \|K^k(f)-\mu(f)\|_{\infty, \mu } \leq
  \|dK^k(f)\| \leq C \rho^k \|df\| \, ,
$$
so that (\ref{bis}) is satisfied. On the other hand, applying
Lemma 1 in Dedecker and Prieur (2007), we have that, for any $l >k
\geq 0$,
\begin{equation}\label{2points}
  \|{\mathbf E}(X_kX_l|Y_0)-{\mathbf E}(X_kX_l)\|_\infty \leq
  C(1+C)\rho^k \|df\|^2 \, ,
\end{equation}
so that (\ref{mix}) holds.  This completes the proof of Corollary
\ref{proexp} when $f$ is $BV$.   \quad $\diamond $

\subsection{Proof of Proposition \ref{circle}}
$\quad \;$ To prove Proposition \ref{circle}, it suffices to prove
that the sequence $X_i= f(\xi_i)-m(f)$ satisfies the conditions
(\ref{bis}) and (\ref{mix}) of Corollary \ref{cor1}. Let
$\|\cdot\|_{\infty, m}$ be the essential supremum norm with respect
to $m$.

Note that $\|{\mathbf{E}}(X_n | \xi_0)\|_\infty= \|K^n(f)- m(f)
\|_{\infty, m}$, and that
\begin{equation*}
K^n(f)(x)-m(f)= \sum_{k \in {\mathbf{Z}}^*} \cos^n(2\pi k a) \hat f(k) \exp
(2i\pi k x) \, .
\end{equation*}
By assumption, there exists $C>0$ such that $\sup_{k \neq 0}
|k|^{1+ \varepsilon}|\hat f (k)| \leq C$. Hence
\begin{equation}  \label{40}
\sum_{n >0} \frac{\|K^n(f)- m(f) \|_{\infty, m}}{\sqrt n} \leq
C\sum_{k \in {\mathbf{Z}}^*} |k|^{-1 - \varepsilon} \sum_{n>0}
\frac{|\cos(2\pi k a)|^n} {\sqrt n} \, .
\end{equation}
Here, note that there exists a positive constant $K$ such that, for
any $0<a<1$, we have $\sum_{n>0} n^{-1/2} a^n \leq K a (1-a)^{-1/2}$
(to see this, it suffices to compare the sum with the integral of
the function $h(x)=x^{-1/2} a^x$). Consequently, we infer from
(\ref{40}) that
\begin{eqnarray}  \label{borne}
\sum_{n >0} \frac{\|K^n(f)- m(f) \|_{\infty, m}}{\sqrt n} &\leq&
CK\sum_{k \in {\mathbf{Z}}^*} \frac{1}{|k|^{1 + \varepsilon}
\sqrt{1-|\cos(2\pi k a)|}}
\notag \\
&\leq & CK \sum_{k \in {\mathbf{Z}}^*}\frac{1}{|k|^{1 + \varepsilon} d(2ka,
{\mathbf{Z}})} \, ,
\end{eqnarray}
the last inequality being true because $(1-|\cos(\pi u)|) \geq \pi
(d(u, {\mathbf{Z}}))^2$. Since $a$ is badly approximable by
rationals, then so is $2a$. Hence, arguing as in the proof of
Lemma 5.1 in Dedecker and Rio (2006), we infer that for any
positive $\eta$ there exists a constant $D$ such that
\begin{equation*}
\sum_{k=2^N}^{2^{N+1}-1} \frac{1}{d(2 k a, {\mathbf{Z}})} \leq D
2^{(N+2)(1+\eta)} N \, .
\end{equation*}
Applying this result with $\eta= \varepsilon/2$, we infer from (\ref{borne})
that
\begin{equation*}
\sum_{n >0} \frac{\|K^n(f)- m(f) \|_{\infty, m}}{\sqrt n} \leq 2CKD
\sum_{N \geq 0} 2^{(N+2)(1+ \varepsilon/2)} N \max_{2^N \leq k \leq
2^{N+1}} k^{-1-\varepsilon} < \infty \, ,
\end{equation*}
so that the condition (\ref{bis}) of Corollary \ref{cor1} is
satisfied. The condition (\ref{mix}) of Corollary \ref{cor1} follows
from the inequality (5.18) in Dedecker and Rio (2006).
\quad $\diamond $

\subsection{Appendix}
$\quad \;$
This section collects some technical lemmas.

The proof of the following lemma is left to the reader since it uses
the same arguments as in the proof of Proposition 2.5 in Peligrad and
Utev (2005) by replacing the $\mathbf{L}_{2}$ norm by the $\mathbf{L}_{\infty}$ norm.

\begin{lemma}
\label{techlemma} Under condition (\ref{mw}),
\begin{eqnarray*}
\frac{\Vert \mathbf{E}(S_{m}|\mathcal{F}_{0})\Vert _{\infty }}{\sqrt{m}}
\rightarrow 0 \quad \text{and} \quad
\frac{1}{\sqrt{m}}\sum_{j=1}^{\infty }\frac{\Vert
\mathbf{E}(S_{mj}|\mathcal{F}_{0})\Vert _{\infty
}}{j^{3/2}}\rightarrow 0 \quad \text{as $m\rightarrow \infty $.}
\end{eqnarray*}
\end{lemma}

The following lemma gives a simple fact about convergence.

\begin{lemma}
\label{convergence} Let $A_{j}(x,n,m)$, $j=1,\ldots,J$, $x>0$, be
real valued functions such that for each $j,n,m$ the function
$A_{j}(x,n,m)$ is non-increasing in $x>0$ and assume that, for any
$x>0$,
\begin{eqnarray*}
\limsup_{m\rightarrow\infty}\limsup_{n\rightarrow\infty}A_{j}(x,n,m)=-\infty
\, .
\end{eqnarray*}
Then for any $u_{n}\rightarrow\infty$, there exists
$m_{n}\rightarrow\infty$ such that $m_{n}\leq u_{n}$ and, for any
$x>0$ and $j=1,\ldots,J$,
\begin{eqnarray*}
\limsup_{n\rightarrow\infty}A_{j}(x,n,m_{n})=-\infty
\end{eqnarray*}
\end{lemma}
{\bf Proof}. First, we observe that by considering the function
\begin{equation*}
A(x,n,m)=\max_{1\leq j\leq J}A_{j}(x,n,m)\;,\;
\end{equation*}
the lemma reduces to the case $J=1$.

Construct a strictly increasing positive integer sequences
$\psi_{k}$ and $n_{k}$ such that for all $n\geq n_{k}$,
\begin{align*}
A(1/k,n,\psi_{k})\leq-k \, .
\end{align*}
Let $g(n)=k$ for $n_{k}<n\leq n_{k+1}$ starting with $k=1$ and $g(n)=1$ for $
n\leq n_{1}$. Then, $g(n)$ is non-decreasing, $g(n)\rightarrow
\infty $ and for all $n>n_{1}$ such that $n_{k}<n\leq n_{k+1}$ (and
so $g(n)=k$).
\begin{equation*}
n_{g(n)}=n_{k}<n
\end{equation*}
Now, let $G(n)$ be a positive integer sequence such that $G(n)\leq g(n)$ and
$G(n)\rightarrow \infty $. Then,
\begin{equation*}
n_{G(n)}\leq n_{g(n)}=n_{k}<n
\end{equation*}
Hence, there exists $G(n)$ such that
\begin{equation*}
\psi _{G(n)}\leq u_{n}\;,\;n_{G(n)}\leq n\quad \mbox{and}\quad
G(n)\rightarrow \infty \;.
\end{equation*}
Finally, let $m_{n}=\psi _{G(n)}$. Then, obviously
\begin{equation*}
m_{n}\leq u_{n}\quad \mbox{and}\quad m_{n}\rightarrow \infty \;.
\end{equation*}
On the other hand, for any $x>0$ and $n$ such that $x\geq 1/G(n)$,
since $A(x,n,m)$ is non-increasing in $x$, we have
\begin{align*}
A(x,n,m_{n})& \leq A(1/G(n),n,m_{n})=A(1/G(n),n,\psi _{G(n)}) \\
& \leq -G(n)\rightarrow -\infty
\end{align*}
which proves the lemma. \quad $\diamond $

\end{document}